\documentclass[11pt,reqno]{amsart}
\usepackage{amsmath}
\usepackage{amsthm}
\usepackage{graphicx}
\usepackage{mathtools}
\usepackage{amsfonts}
\usepackage{amssymb}
\usepackage{hyperref}
\usepackage{bm}
\usepackage[margin=1 in]{geometry}
\usepackage{enumerate}
\usepackage[normalem]{ulem}
\usepackage{tikz}
\usepackage{tikz-cd}
\usepackage{xcolor}
\usetikzlibrary{matrix,arrows,positioning,automata}
 \usepackage{cancel}

\usepackage{dsfont}

\newcommand{\mdiv}{{\rm div}}
\newcommand{\R}{\mathbb{R}}

\newcommand{\E}{\mathbb{E}}

\renewcommand{\emptyset}{\varnothing}

\newcommand{\cL}{\mathcal{L}}

\newtheorem{theorem}{Theorem}[section]
\newtheorem{corollary}{Corollary}[section]
\newtheorem{lemma}{Lemma}[section]
\newtheorem{proposition}{Proposition}[section]

\theoremstyle{remark}
\newtheorem{remark}{Remark}

\theoremstyle{definition}
\newtheorem{definition}{Definition}[section]

\newtheorem{assumption}{Assumption}

\title{Mean-field limit for stochastic control problems under state constraint}
\author[S. Daudin]{Samuel Daudin 
\address{(S. Daudin) Universit\'e C\^ote d'Azur, CNRS, Laboratoire J.A. Dieudonn\'e, 06108 Nice, France
}\email{samuel.daudin@univ-cotedazur.fr
}}

\begin{document}


\maketitle

\begin{abstract}
We study the convergence problem of mean-field control theory in the presence of state constraints and non-degenerate idiosyncratic noise. Our main result is the convergence of the value functions associated to stochastic control problems for many interacting particles subject to symmetric, almost-sure constraints toward the value function of a control problem of mean-field type, set on the space of probability measures. The key step of the proof is to show that admissible controls for the limit problem can be turned into admissible controls for the $N$-particle problem up to a correction which vanishes as the number of particles increases. The rest of the proof relies on compactness methods. We also provide optimality conditions for the mean-field problem and discuss the regularity of the optimal controls. Finally we present some applications and connections with large deviations for weakly interacting particle systems.
\end{abstract}

\tableofcontents

\section*{Introduction}

The goal of this paper is to address the so-called convergence problem of mean-field control theory in the presence of state constraints and non-degenerate idiosyncratic noise. The pre-limit problem involves a large number of interacting particles subject to symmetric, almost-sure constraints while, in the limit, the constraint acts on the law of a typical particle. Constraints in law arise naturally in applications in economy and finance, as a way to control the risk associated with a given strategy, \cite{Follmer1999, Krokhmal2001, Markowitz1952}. They have generated some recent attention in the stochastic control community \cite{Bouchard2010,Chow2020,Frankowska2018,Frankowska2019, Germain2022, Pfeiffer2020, Pfeiffer2020a} and the development of mean-field game and mean-field control theory \cite{Caines2006,Caines2007,Cardaliaguet2019a, Carmona2018a, Carmona2018b,Lasry2006,Lasry2006a,Lasry2007,Lions20062020} provides new insights and techniques to address these problems. Notably, given the nature of the constraint, we are naturally led to consider control problems on the space of probability measures which require appropriate tools, \cite{Bonnet2019, Bonnet2021, Carmona2018a}. Constraints in law also arise in continuous descriptions of controlled particle systems \cite{Cardaliaguet2016,DiMarino2016b, Meszaros2015,Meszaros2018,Meszaros2016a, Seguret2021}. While the convergence problem in mean-field control is well understood \cite{Budhiraja2012, Cardaliaguet2022, Cardaliaguet2022a, Daudin2023a, Fischer2016, Germain2021, Lacker2017}  the goal of this paper is to investigate the validity of this approximation in the presence of constraints.

More precisely we investigate the connection between the following two control problems. The first one involves a large number $N\geq 1$ of interacting particles. Its value function is given by 
\begin{equation}
\mathcal{V}^N \bigl(t_0,\mathbf{x}^N_0 \bigr) = \inf_{(\alpha_t^{i,N})_{1\leq i \leq N}} \E \left[ \int_{t_0}^T \frac{1}{N} \sum_{i=1}^N L \bigl(X_t^{i,N}, \alpha_t^{i,N} \bigr) dt + \int_{t_0}^T \mathcal{F} \bigl( \hat{\mu}_t^{N} \bigr) dt +\mathcal{G} \bigl(\hat{\mu}_T^{N} \bigr)\right]
\label{NPIntro2022}
\tag{NP}
\end{equation}
where $T>0$ is a finite horizon, $t_0 \in [0,T]$ is the initial time and $\mathbf{x}^N_0=(x_0^{1,N},\dots,x_0^{N,N}) \in (\R^d)^N$ denotes the initial position of the particles. The dynamics are given by the stochastic differential equations


\begin{minipage}{0.48 \textwidth}
\begin{equation}
\left \{ 
\begin{array}{ll} 
\displaystyle  dX_t^{i,N} =b \bigl(X_t^{i,N}, \hat{\mu}_t^{N} \bigr)dt +  \displaystyle \alpha_t^{i,N} dt + \sqrt{2}dB_t^{i,N} \\
\bigl(X_{t_0}^{1,N},\dots, X_{t_0}^{N,N} \bigr) =\mathbf{x}^N_0
\end{array}
\right.
\end{equation}
\end{minipage}
\begin{minipage}{0.48 \textwidth}
\begin{equation}
\displaystyle  \hat{\mu}_t^{N} := \frac{1}{N} \sum_{i=1}^N \delta_{X_t^{i,N}} 
 \end{equation}
\end{minipage}
and the infimum is taken over a suitable class of controls $\mathbf{\alpha} =\bigl(\alpha^{1,N}, \dots, \alpha^{N,N} \bigr)$. Importantly, the particles $(X^{1,N},\dots, X^{N,N})$ are subject to the state constraint
$$\Psi^N \bigl( X_t^{1,N} , \dots, X_t^{N,N} \bigr) < 0 \mbox{ for all } t \in [0,T], \quad\mathbb{P}-\mbox{almost-surely}, $$
for some symmetric map $\Psi^N : (\R^d)^N \rightarrow \R$. We will always assume that $\Psi^N$ has the form
$$\displaystyle \Psi^N \bigl( x^{1,N}, \dots, x^{1,N} \bigr) = \Psi \Bigl( \frac{1}{N} \sum_{i=1}^N \delta_{x^{i,N}} \Bigr), $$  
for some functionnal $\Psi$ defined over $\mathcal{P}(\R^d)$ the set of Borel probability measures over $\R^d$. Therefore the constraint reads
$$\Psi \bigl( \hat{\mu}_t^{N} \bigr ) < 0 \mbox{ for all } t \in [0,T], \quad \mathbb{P}-\mbox{almost-surely}. $$
Above, the cost function $L : \R^d \times \R^d \rightarrow \R$ is supposed to be convex with quadratic growth in the second variable and smooth. The non-linear drift $b : \R^d \times \mathcal{P}_1(\R^d) \rightarrow \R^d$ is (at least) Lipschitz continuous and bounded. The mean-field costs $\mathcal{F}, \mathcal{G} : \mathcal{P}_1(\R^d) \rightarrow \R$ as well as the constraint $\Psi : \mathcal{P}_1(\R^d) \rightarrow \R$ are (at least) Lipschitz continuous functions over the space of Borel probability measures with finite first order moment $\mathcal{P}_1(\R^d)$ (precise assumptions will be given later in Section \ref{sec: Assumptions2023} ). For each $N \geq 1$, the $(B^{i,N})_{ 1\leq i \leq N}$ are independent $d$-dimensional Brownian motions.


The second problem of interest in this paper in an optimal control problem for a non-linear Fokker-Planck equation. Its value function is defined for every $(t_0,m_0) \in [0,T] \times \mathcal{P}_2(\R^d)$ such that $\Psi(m_0) \leq 0$ by
\begin{equation} 
 \mathcal{U}(t_0,m_0) = \inf_{(\mu,\alpha)} \int_{t_0}^T \int_{\R^d} L\bigl(x,\alpha(t,x) \bigr)d\mu(t)(x)dt + \int_{t_0}^T \mathcal{F}\bigl(\mu(t)\bigr) dt + \mathcal{G}\bigl(\mu(T)\bigr) 
\label{MFproblemIntro2022}
\tag{mfP}
\end{equation}
where the infimum is taken over the couples $(\mu,\alpha)$ with $\mu \in \mathcal{C}([t_0,T] ,\mathcal{P}_2(\R^d)), \alpha\in L^2_{dt \otimes \mu(t)} ([t_0,T] \times \R^d,\R^d) $ satisfying the non-linear Fokker-Planck equation
\begin{equation}
\left \{
\begin{array}{ll}
\partial_t \mu + \mdiv\bigl((\alpha(t,x) + b(x,\mu(t))\mu \bigr) -\Delta \mu = 0 &\mbox{ in } (t_0,T) \times \R^d,  \\
\mu(t_0)=\mu_0 \in \mathcal{P}_2(\R^d),
\end{array}
\right.
\end{equation}
subject to the state constraint
$$\Psi(\mu(t)) \leq 0, \quad \forall t \in [t_0,T]. $$

The latter problem was analyzed by the author in \cite{Daudin2021}. It is proved there, under appropriate conditions on the data involving qualification conditions for the constraint, that optimal controls are bounded and Lipschitz continuous in space, at least for any initial position $\mu_0 \in \mathcal{P}_2(\R^d)$ such that $\Psi(\mu_0) <0$. Problem \eqref{NPIntro2022} which defines $\mathcal{V}^N$ is, however, very different in nature. Indeed the constraint has to be satisfied almost-surely while the noises driving the dynamics of the particles are non-degenerate. This type of constraint leads to new difficulties. Indeed, to dominate the effect of the diffusion, controls cannot remain bounded and the value function associated to this problem blows-up near the boundary. 


Without constraint, the connection between Problem \eqref{MFproblemIntro2022} and Problems \eqref{NPIntro2022} is by now well understood. Under more general structure conditions, Lacker proved in \cite{Lacker2017} that the law of the empirical measures of weak solutions to the $N$-particle system converges to probability measures supported on the set of optimal solutions to the mean-field problem and therefore convergence of the value functions hold. Taking advantage of the regularizing effect of the diffusion and uniform in $N$ Lipschitz and semi-concavity estimates for the value functions of the $N$-particles system, it was shown in \cite{Cardaliaguet2022} that convergence actually holds with a rate. In the same setting, Cardaliaguet and Souganidis later proved in \cite{Cardaliaguet2022a} a propagation of chaos around “stable" solutions of the mean-field problem. We mention that, under convexity assumptions on the mean-field costs $\mathcal{F}$ and $\mathcal{G}$ it is shown in \cite{Cardaliaguet2019a} that the value function associated to the mean-field control problem is a smooth (enough) function in the Wasserstein space. In this setting it is not difficult to prove that the convergence of the value functions holds with an optimal rate and we have quantitative propagation estimates for the optimal trajectories to the $N$-particles system toward the solution to the mean-field control problem. Finally, the recent contribution \cite{Daudin2023a} by the author, together with Delarue and Jackson, provides optimal rates of convergence under appropriate regularity conditions on the data but without assuming that the value function $\mathcal{U}$ is differentiable.

We also mention that recent progresses were made in order to characterize the value function for the mean-field problem, in the general situation where it is not expected to be smooth. Similarly to the finite dimensional case, we expect the value function to be the unique viscosity solution (in some sense) to the dynamic programming equation. Different approaches have been taken in \cite{Burzoni2020, Cecchin2022, Conforti2021, Cosso2021}. The most general result, so far, being \cite{Cosso2021}, where the authors rely on the approximation of the mean-field control problem by control problems for finite numbers of interacting particles. An analog characterization of the value function in the presence of state constraints is still an open question.

Stochastic control problems with state constraints and non-degenerate diffusions were addressed in the seminal work \cite{Lasry1989} of Lasry and Lions. They showed that the blow-up behavior of the value function is directly related to the growth of the Hamiltonian and provided rates of divergence. This problem was later revisited by Leonori and Porretta in \cite{Leonori2007} where the authors also find the rate of divergence of the optimal controls. Another approach consists in requiring the volatility to degenerate at the boundary in order to find admissible controls with  a finite cost, see for instance \cite{Bouchard2009}. In this case, Dynamic Programming leads to constrained viscosity solutions in the sense of Soner \cite{Soner1986} to the corresponding Hamilton Jacobi Bellman equation, see \cite{Katsoulakis1994a}.

Without constraint, Problem \eqref{MFproblemIntro2022} has been widely studied in the literature, mainly for its connection with potential mean-field games. Indeed optimality conditions for Problem \eqref{MFproblemIntro2022} (without constraint) lead to the celebrated mean-field game system of PDEs introduced by Lasry and Lions in \cite{Lasry2007}. We refer to \cite{Briani2018,Cardaliaguet2019a,Lasry2007} for various properties of the pde system in the general case and to \cite{Daudin2021} in the presence of state constraints.

In the context of mean-field control, state constraints have been primarily studied at the level of the limit problem in order to derive versions of the Pontryagin maximum principle. This is achieved in \cite{Bonnet2019, Bonnet2021} for first order problem (namely without diffusion). In  \cite{Frankowska2019} the authors provide first and second order conditions for optimality for stochastic control problems with expectation constraints, which corresponds to Problem \eqref{MFproblemIntro2022}. While much effort was made to understand optimal control problems over the space of probability measures under state constraints, less is known about the connection between the limit problem and the associated $N$-particle problem. This is a particularly challenging question when the dynamics of the particles are stochastic.

In this paper we prove the convergence of the value functions for the problems with almost-sure constraints toward the value function for the mean-field problem. Similarly to \cite{Cardaliaguet2022} we proceed in two steps. On the one hand we prove that 
$$  \mathcal{U}(t_0,\mu_0) \leq \liminf_{N \rightarrow +\infty} \mathcal{V}^N \bigl(t_0,x_0^{1,N},\dots,x_0^{N,N} \bigr),$$
whenever $\hat{\mu}_0^N$ converges to $\mu_0$ in $\mathcal{P}_2(\R^d)$. This boils down to finding weak limit points of sequences of nearly optimal weak solutions to the $N$-particle problem. Once we know that $\mathcal{V}^N$ is bounded independently from $N$, this follows from the line of arguments of \cite{Lacker2017} for problems without constraint. 
On the other hand, proving that 
$$\limsup_{N \rightarrow +\infty} \mathcal{V}^N\bigl(t_0,x_0^{1,N},\dots,x_0^{N,N}\bigr) \leq \mathcal{U}(t_0,\mu_0) $$
requires more care. Indeed an admissible control for the mean-field problem is, in general, not admissible for the particle system because of the almost-sure constraint.


Our strategy can be described as follows. Given an admissible control $\alpha$ for the mean-field control problem, we consider the particle system starting from an initial position $\bigl(X_0^{1,N},\dots,X_N^{N,N}\bigr) =\bigl(x_0^{1,N},\dots, x_0^{N,N}\bigr) $,
$$ X_t^{i,N} = x_0^{i,N} + \int_{t_0}^{t \wedge \tau_N} \alpha\bigl(s,X_s^{i,N}\bigr) dt  +\int_{t_0}^{t \wedge \tau_N} b\bigl(X_s^{i,N},\hat{\mu}_s^{N}\bigr)ds +\int_{t \wedge \tau_N }^t \beta_t^{i,N} dt + \sqrt{2}\int_{t_0}^{t \wedge \tau_N} dB_t^{i,N} $$
where $\tau_N := \inf \bigl \{t \geq t_0, \Psi( \hat{\mu}_t^{N} ) \geq -\frac{\delta }{2} \bigr \}$ and $\beta_t^{i,N}$ is a feedback control designed so that, $\mathbb{P}$-almost-surely
$$ \frac{1}{N} \sum_{i=1}^N \left| X_t^{i,N} - X_{\tau_N}^{i,N} \right|^2 \leq r^2 \hspace{30pt} \forall t \geq \tau_N, $$
where $r$ is a small radius depending on $\delta$ which guarantees that, $\mathbb{P}$-almost-surely,
$$ \Psi( \hat{\mu}_t^{N} ) < 0, \hspace{30pt} \forall t\in [t_0,T]. $$
If $\alpha$ is bounded, Lipschitz continuous and taken so that the corresponding solution $\mu$ to 
\begin{equation}
\left \{
\begin{array}{ll}
\partial_t \mu + \mdiv\bigl((\alpha(t,x) + b(x,\mu(t))\mu \bigr) -\Delta \mu = 0 & \mbox{ in } (t_0,T) \times \R^d, \\
\mu(t_0)=\mu_0
\end{array}
\right.
\end{equation}
satisfies $\Psi(\mu(t)) \leq -\delta$, for all $t \in [t_0,T]$, for some $\delta >0$, we expect a strong convergence of $\hat{\mu}_t^{N}$  toward $\mu(t)$ for $t \in [t_0, \tau_N]$ and therefore $\tau_N \wedge T$ must converge to $T$. The key step is to build $\bigl(\beta_t^{i,N}\bigr)_{1 \leq i \leq N}$ so that its contribution to the cost for the $N$-particle problem, vanishes as $N \rightarrow +\infty$. We are able to do so only if $\Psi(\mu_0)<0$.
We also need to prove that it is enough to approximate admissible candidates $(\alpha, \mu)$ such that $\alpha$ is bounded, Lipschitz continuous with respect to the space variable and $\Psi(\mu(t)) \leq -\delta$ for all $t \in [t_0,T]$, for some $\delta >0$. 

Overall, our main result, Theorem \ref{MainTheoremMFconstraint2022}, states that, under Assumption \eqref{AssumptionMeanFieldLimitSept2022} (introduced in Section \ref{sec: Assumptions2023}), we have 
$$ \lim_{N \rightarrow +\infty} \mathcal{V}^N \bigl(t_0,x_0^{1,N},\dots,x_0^{N,N} \bigl) = \mathcal{U}(t_0,\mu_0), $$
whenever $\frac{1}{N}\sum_{i=1}^N \delta_{x_0^{i,N}}$ converges in $\mathcal{P}_2(\R^d)$ toward some $\mu_0$ such that $\Psi(\mu_0)<0$. We also prove that accumulation points of weak solutions (introduced in Subsection \ref{sec: The N-state problem}) to the $N$-particle problems are supported on the set of solutions of the limiting problem and therefore we have a proper convergence result for the optimal solutions when the limit problem admits a unique solution. \\

\textbf{Connection with large deviations for weakly interacting diffusions.} Our result is closely related to the large deviation principle for weakly interacting particle systems, see \cite{Budhiraja2000,Dawson1987,Fischer2014}. Indeed, consider the (uncontrolled) particle system 
\begin{equation*}
\left \{
\begin{array}{ll}
\displaystyle X_t^{i,N} = x^{i,N} + \int_{0}^t b \bigl(X_s^{i,N}, \hat{\mu}_s^{N}\bigr) ds +\sqrt{2}B_t^{i,N}, &  t \in [0,T], \quad i \in \left \{1,\dots,N \right \}, \\
\displaystyle \hat{\mu}_s^{N} = \frac{1}{N} \sum_{i=1}^N \delta_{X_s^{i,N}}, 
\end{array}
\right.
\end{equation*}
where $(B_t^{1,N}, \dots, B_t^{N,N})$ are $N$ independent $d$-dimensional standard Brownian motions supported on some probability space, and let $v^N(t, \mathbf{x}^N)$ be  the probability that the particles initialized at $(0, \mathbf{x}^N)$ stay strictly inside the constraint, at least up to time $t$. That is 
\begin{align*}
v^N(t,\mathbf{x}^N) &:= \mathbb{P} \bigl( \forall s \in [0,t], \Psi (\hat{\mu}_s^N) <0 \bigr).
\end{align*}
Under appropriate assumptions on $b$ and $\Psi$, the map $ (t,\mathbf{x}^N) \mapsto -\frac{2}{N} \log v^N(T-t,\mathbf{x}^N), $ solves
the dynamic programming equation for Problem \eqref{NPIntro2022}  when $\mathcal{F} = \mathcal{G}=0$ and $L(x,q) = \frac{1}{2}|q|^2$. As a consequence we can deduce the exponential decay of $v^{N}$ by looking at the limit of $\mathcal{V}^N(t, \mathbf{x}^N)$ when $\displaystyle \frac{1}{N}\sum_{i=1}^N \delta_{x^{i,N}} \rightarrow \mu_0$ in $\mathcal{P}_2(\R^d)$. In Section \ref{sec: Application to Large Deviations Sept2022} we discuss this rigorously. Notice that this method to obtain estimates on the probability $v^N$ by making a logarithmic transformation and studying the stochastic control problem corresponding to the resulting Hamilton-Jacobi-Bellman equation is reminiscent of \cite{Fleming1977}.\\

\paragraph{\textbf{Organization of the paper}} The rest of  the paper is organized as follows. In Section \ref{sec: Assumptions and main results 11MAI2023} we give our working assumptions, the precise formulation of the problems considered and our main results. In Section \ref{sec: Properties of the mean-field problem 11 mai 2023} we give properties of the mean-field problem. In particular, we give optimality conditions in the form of a mean-field game system of pdes. We also discuss essential stability properties of the men-field problem with respect to small perturbations of the constraint. In the next section, Section \ref{sec: Mean field limit 2022}, we prove our main convergence result. In Section \ref{sec: Application to Large Deviations Sept2022} we present some applications to theory of large deviations for weakly interacting particles. Finally we postpone to Appendix \ref{subsec:OptimalityConditions} the proof of the optimality conditions of Section \ref{sec: Properties of the mean-field problem 11 mai 2023} and to Appendix \ref{sec: Concentration limit} a useful lemma about the convergence of weakly interacting particle systems (without control). \\


\paragraph{\textbf{Notation}} The Wasserstein space of Borel probability measures over $\R^d$ with finite moment of order $r \geq 1$ is denoted by $\mathcal{P}_r(\R^d)$. It is endowed with the $r$-Wasserstein distance $d_r$. 

For $n \geq 1$ we denote by $E_n$ the subspace of $\mathcal{C}^n(\R^d)$ consisting of functions $u$ such that 
$$ \|u\|_n := \sup_{x \in \R^d} \frac{|u(x)|}{1+|x|} + \sum_{k=1}^n \sup_{x \in \R^d} \bigl|D^k u(x) \bigr| <+ \infty. $$
Similarly we define $E_{n+\alpha}$ for $n \geq 1$ and $\alpha \in (0,1)$ to be the subset of $E_n$ consisting of functions $u$ satisfying
$$ \|u \|_{n +\alpha} := \|u \|_n + \sup_{x \neq y} \frac{|D^n u(x) - D^n u(y) |}{|x-y|^{\alpha}} < +\infty. $$ 

Finally we will use the heat kernel $P_t$ associated to $-\Delta$ defined, when it makes sense, by
$$ P_t f (x) := \int_{\R^d} \frac{1}{(4 \pi t)^{d/2}}e^{-\frac{|x-y|^2}{4t}}f(y)dy. $$

\section{Assumptions and statement of the main results}

\label{sec: Assumptions and main results 11MAI2023}

\subsection{Assumptions}

\label{sec: Assumptions2023}

We first give the assumptions satisfied by $L$, $\mathcal{F}$, $\mathcal{G}$ and $\Psi$. They involve an integer $n\geq 3$. For $U=\mathcal{F}, \mathcal{G}, \Psi$, the map $U: \mathcal{P}_2(\R^d) \rightarrow \R^d$ satisfies 
\begin{equation}
U \mbox{ is a bounded from below, } \mathcal{C}^1\mbox{ map and } \frac{\delta U}{\delta m} \mbox{  belongs to }\mathcal{C}(\mathcal{P}_2(\R^d),E_{n+\alpha}),
\tag{Ureg}
\label{Ureg}
\end{equation}
We recall that  a map $U : \mathcal{P}_2(\R^d) \rightarrow \R^m$ is $\mathcal{C}^1$ if there exists a jointly continuous map $\displaystyle \frac{\delta U}{\delta m} : \mathcal{P}_2(\R^d) \times \R^d \rightarrow \R^m$ such that, for any bounded subset $\mathcal{K} \subset \mathcal{P}_2(\R^d)$, $\displaystyle x \rightarrow \frac{\delta U}{\delta m}(m,x)$ has at most quadratic growth in $x$ uniformly in $m \in \mathcal{K}$ and such that, for all $m, m' \in \mathcal{P}_2(\R^d)$,
$$U(m')-U(m)= \int_0^1 \int_{\R^d} \frac{\delta U}{\delta m} \bigl((1-h)m+hm',x \bigr)d(m'-m)(x)dh. $$
The function $\displaystyle \frac{\delta U}{\delta m}$ is defined up to an additive constant and we adopt the normalization convention $\displaystyle \int_{\R^d} \frac{\delta U}{\delta m}(m,x)dm(x) = 0.$ We refer to the monographs \cite{Cardaliaguet2019a,Carmona2018a} for a discussion about the notion(s) of derivatives in the space of probability measures.

The Lagrangian $L$ verifies $L(x,q) = \sup_{p\in \R^d} \left \{ -p.q-H(x,p) \right \}$ for all $(x,q) \in \R^d \times \R^d$ where $H$, the Hamiltonian, satisfies the following conditions.
\begin{equation}
\left \{
\begin{array}{ll}
\displaystyle \mbox{$H$ belongs to $\mathcal{C}^n (\R^d \times \R^d)$.} \\
\displaystyle \mbox{$H$ and its derivatives are bounded on sets of the form $\R^d \times B(0,R)$ for all $R >0$.} \\
\displaystyle \mbox{For some $C_0>0$, for all $(x,p) \in \R^d \times \R^d$}, \\
\displaystyle \hspace{150pt} |D_xH(x,p)| \leq C_{0} (1 +|p|). \\
\displaystyle \mbox{For some $\mu >0$ and all $(x,p) \in \R^d \times \R^d$,} \\
\hspace{150pt}  \frac{1}{\mu} I_d \leq D^2_{pp}H(x,p) \leq \mu I_d.
\end{array}
\right.
\tag{AH}
\label{AHallinone}
\end{equation}
The non-linear drift $b$ is assumed to satisfy
\begin{equation}
\begin{array}{cc}
b: \R^d \times \mathcal{P}_2(\R^d) \rightarrow \R \mbox{ is bounded, Lipschitz continuous and admits a linear derivative } \\
\displaystyle \frac{\delta b}{\delta m} \in \mathcal{C}_b \bigl(\R^d \times \mathcal{P}_2(\R^d), E_n(\R^d) \bigr).
\end{array}
\tag{Ab}
\label{Ab}
\end{equation}
For the constraint, we also assume that 
\begin{equation}
\Psi \mbox{ is convex}, 
\label{APsiConv}
\tag{APsiConv}
\end{equation}
Finally we also assume that:
\begin{equation}
\mbox{ There is at least one } \mu \in \mathcal{P}_2(\R^d) \mbox{  such that } \Psi(\mu)<0.
\label{APsiInside2023}
\tag{APsiInside}
\end{equation}
For convenience we put all of the above assumptions into 
\begin{assumption}
Assume that \eqref{Ureg} holds for $\mathcal{F},\mathcal{G}$ and $\Psi$, \eqref{AHallinone} holds for $H$, \eqref{Ab} holds for $b$ and \eqref{APsiConv}, \eqref{APsiInside2023} hold for $\Psi$.
\label{AssumptionMeanFieldLimitSept2022}
\end{assumption}

A typical example of functions satisfying the condition \eqref{Ureg} is the class of cylindrical functions of the form 
$$ U(m) = F \left( \int_{\R^d} f_1(x)dm(x), \dots , \int_{\R^d} f_k(x)dm(x) \right), $$
where $F$ and the $f_i$, $1 \leq i \leq k$ are smooth with bounded derivatives.  Assumption \eqref{Ureg} also implies that $(m,x) \rightarrow D_m U(m,x)$ is uniformly bounded in $\mathcal{P}_2(\R^d) \times \R^d$ and therefore, a simple application of Kantorovitch-Rubinstein duality for $d_1$ proves that $U$ is Lipschitz continuous with respect to this distance. 
Finally, $\Psi(m) := \int_{\R^d} \psi(x)dm(x)$ satisfies Assumptions \eqref{APsiConv} and \eqref{APsiInside2023} as soon as $\psi \in E_n$ (not necessarily convex) and there is $x_0 \in \R^d$ such that $\psi(x_0)<0$. In this “linear" case we recover the control problem with expectation constraints, see \cite{Bouchard2010, Chow2020,Frankowska2018, Frankowska2019, Pfeiffer2020a}.

\subsection{The problem with almost-sure constraint}

\label{sec: The N-state problem}

\hspace{100pt}

\textbf{Strong formulation}. Throughout this section we fix some $t_0 \in [0,T]$ and $\mu_0 \in \mathcal{P}_2(\R^d)$ such that $\Psi(\mu_0) <0$. In its strong formulation, the $N$-state control problem is described as follows. We fix a probability space $(\Omega, \mathcal{F},\mathbb{P})$ endowed with $N$ independent standard Brownian motions $(B_t^{i,N})_{i=1,...,N}$. We also fix some initial positions $\mathbf{x}_0 = \bigl(x_0^{1,N},\dots,x_0^{N,N} \bigr) \in (\R^d)^N$ such that $\Psi \left( \frac{1}{N}\sum_{i=1}^N \delta_{x_0^{i,N}} \right) <0$ for all $N \geq 1$.

The controller's problem is to minimize over controls $(\alpha_t^{i,N})_{i=1,\dots,N} \in L^2([0,T] \times \Omega, (\R^d)^N)$ adapted to the filtration generated by the Brownian motions
$$J^N \bigl(t_0, \mathbf{x}_0;(\alpha_t^{i,N})_{1 \leq i \leq N} \bigr):= \mathbb{E} \left[ \int_{t_0}^T  \left( \frac{1}{N} \sum_{i=1}^N L \bigl(X_t^{i,N},\alpha_t^{i,N}\bigr) + \mathcal{F} \bigl(\hat{\mu}_t^{N} \bigr)\right)dt + \mathcal{G} \bigl(\hat{\mu}_T^{N} \bigr) \right] $$
under the dynamics
\begin{equation*}
\left \{
\begin{array}{ll}
\displaystyle X_t^{i,N} = x_0^{i,N} +\int_{t_0}^t b \bigl(X_s^{i,N}, \hat{\mu}_s^{N} \bigr)ds + \int_{t_0}^t \alpha_s^{i,N} ds + \sqrt{2} \bigl(B_t^{i,N} - B_{t_0}^{i,N} \bigr), \hspace{30pt} t \geq t_0, \\
\displaystyle \hat{\mu}_t^{N}:= \frac{1}{N}\sum_{i=1}^N \delta_{X_t^{i,N}}
\end{array}
\right.
\end{equation*}
and the constraint
\begin{equation} 
\Psi(\hat{\mu}_t^{N}) < 0, \hspace{20pt} \mbox{for all  } t \in[t_0,T], \hspace{20pt}  \mathbb{P}-\mbox{almost-surely}. 
\label{ConstraintNparticle10Mai}
\end{equation} 
We will also denote the constraint for the $N$-particle problem by
$$\Omega_N:= \left \{ \mathbf{x}^N \in \bigl(\R^d \bigr)^N, \Psi \Bigl( \frac{1}{N}\sum_{i=1}^N\delta_{x^{i,N}} \Bigr)<0 \right \}$$
and, in this case, \eqref{ConstraintNparticle10Mai} reads
$$ \bigl(X_t^{1,N}, \dots, X_t^{N,N} \bigr) \in \Omega_N, \hspace{20pt} \mbox{for all  } t \in[t_0,T], \hspace{20pt}  \mathbb{P}-\mbox{almost-surely}. $$
We denote by $\mathcal{V}^N(t_0,\mathbf{x}_0)$ the value of the above problem. Under Assumption \eqref{AssumptionMeanFieldLimitSept2022} we expect, by dynamic programming or a verification argument ---see \cite{Lasry1989,Leonori2007}--- , that $\mathcal{V}^N$ satisfies the Hamilton-Jacobi-Bellman equation 
\begin{equation}
\left \{
\begin{array}{ll}
\displaystyle -\partial_t \mathcal{V}^N -\sum_{i=1}^N b^{i,N}(\mathbf{x}^N).D_{x^{i,N}}\mathcal{V}^N\\
\displaystyle \hspace{60pt}+ \frac{1}{N}\sum_{i=1}^N H\bigl( x^{i,N},ND_{x^{i,N}}\mathcal{V}^N \bigr)  -\sum_{i=1}^N\Delta_{x^{i,N}} \mathcal{V}^N = \mathcal{F}^N(\mathbf{x}^N), & \mbox{ in } (0,T) \times \Omega_N \\
\displaystyle \mathcal{V}^N(t,\mathbf{x}^N) = +\infty, & \mbox{ in } [0,T] \times \partial \Omega_N \\
\displaystyle \mathcal{V}^N(T,\mathbf{x}^N) =\mathcal{G}^N(\mathbf{x}^N) & \mbox{in } \Omega_N,
\end{array}
\right.
\label{EquationforVN11Mai}
\end{equation}
where $\mathcal{F}^N: \bigl( \R^d \bigr)^N \rightarrow \R$ and $\mathcal{G}^N: \bigl( \R^d \bigr)^N \rightarrow \R$ are defined by
$$ \mathcal{F}^N \bigl( \textbf{x}^N \bigr) = \mathcal{F} \Bigl( \frac{1}{N} \sum_{i=1}^N \delta_{x^{i,N}} \Bigr), \quad \mathcal{G}^N \bigl( \textbf{x}^N \bigr) = \mathcal{G} \Bigl( \frac{1}{N} \sum_{i=1}^N \delta_{x^{i,N}} \Bigr), \quad \quad \textbf{x}^N \in\bigl( \R^d \bigr)^N.$$   
In the special case $H(x,p) = 1/2 |p|^2$ and under additional assumptions on $\Psi$ we actually prove that $\mathcal{V}^N$ belongs to $\mathcal{C}^{1,2}\bigl( [0,T) \times \Omega_N \bigr)$ and satisfies \eqref{EquationforVN11Mai} --- see the comment after Proposition \ref{ControlRepresentation} in Section \ref{sec: Application to Large Deviations Sept2022}.

\begin{remark}
Notice that it could very well happen that $\Omega_N = \emptyset$ for small values of $N$. However we neglect this detail since we always assume that there is some $\mu_0 \in \mathcal{P}_2(\R^d)$ such that $\Psi(\mu_0) <0$. By an approximation argument we find that $\Omega_N$ is not empty for $N$ large enough.
\end{remark}

\textbf{Weak formulation.} Let us introduce some notation. We denote by $\mathcal{C}^d:= \mathcal{C}([t_0, T], \R^d)$ the path space. The control space $\mathcal{V}$ is defined as the set of non-negative measures $q$ over $[t_0,T] \times \R^d$ with the Lebesgue measure as time marginal and such that 
$$ \int_{[t_0, T] \times \R^d} |a|^2 dq(t,a) < +\infty. $$

We denote by $(X,\Lambda)$ the canonical process on $(\mathcal{C}^d \times \mathcal{V} )$ and by $(X^{i,N}, \Lambda^{i,N})_{1\leq i \leq N}$ the canonical process on $(\mathcal{C}^d \times \mathcal{V})^N$ and define the empirical measures
$$ \hat{\nu}^N := \frac{1}{N}\sum_{i=1}^N \delta_{(X^{i,N}, \Lambda^{i,N})}, \hspace{30pt} \hat{\mu}^{N}_t := \frac{1}{N}\sum_{i=1}^N \delta_{X_t^{i,N}} = X_t\#\nu^N.$$
We define $\mathcal{R}^N$ as the set of probabilities $P_N \in \mathcal{P}_2((\mathcal{C}^d \times \mathcal{V})^N)$ under which $(X^{i,N}_0)_{i=1,\dots,N}=\mathbf{x}_0^N$, $P_N$-almost-surely,
$$\varphi \bigl(X_t^{1,N},\dots, X_t^{N,N} \bigr) - \sum_{i=1}^N\int_{t_0}^t \int_{\R^d} \mathcal{L}_i^N \varphi \bigl(\hat{\mu}_s^{N},X_s^{1,N}, \dots, X_s^{N,N},a \bigr) d\Lambda^{i,N}_s(a)ds $$
is a martingale under $P_N$, for all smooth, compactly supported $\varphi$ with 
$$ \mathcal{L}_i^N \varphi (\mu,x_1,\dots,x_N,a) := D_{x_i} \varphi(x_1,\dots,x_N).a + D_{x_i} \varphi(x_1,\dots,x_N).b(x_i,\mu)+ \Delta_{x_i}\varphi(x_1,\dots,x_N). $$
The control rule $P_N$ is also assumed to satisfy
$$P_N \Bigl( \Psi( \hat{\mu}^N_{t} ) < 0, \quad \forall t \in [t_0,T] \Bigr) = 1.$$
The $N$-state problem in its weak formulation is therefore to minimize over $P_N \in \mathcal{R}^N$ the cost functional 
$$\mathbb{E}^{P_N} \left [ \int_{t_0}^T \left(  \int_{\R^d}  \frac{1}{N } \sum_{i=1}^N L(X_t^{i,N}, a) d\Lambda^{i,N}_t(a) + \mathcal{F}(\hat{\mu}^N_{t}) \right)dt  + \mathcal{G}( \hat{\mu}_T^{N}) \right]$$
where $\mathbb{E}^{P_N}$ is the expectation under $P_N$. 
We denote by $\mathcal{V}^N_{w}(t_0, \textbf{x}_0)$ the value of this new problem. The next result asserts that the value of the weak formulation is no greater than the value of the strong one.

\begin{lemma}
For all $(t_0, \mathbf{x}^N_0) \in [0,T] \times (\R^d)^N$ such that $\Psi \Bigl( \frac{1}{N}\sum_{i=1}^N \delta_{x_0^{i,N}} \Bigr) < 0 $ it holds
$$\mathcal{V}^N_{w}(t_0, \mathbf{x}^N_0) \leq \mathcal{V}^N(t_0, \mathbf{x}^N_0).$$
\label{lem:weakvsstrongNparticleMai2023}
\end{lemma}

\begin{proof}
It suffices to consider an admissible control $\bigl(\alpha^{i,N}_t \bigr)_{t_0 \leq t \leq T}, 1\leq i \leq N$ for $\mathcal{V}^N \bigl(t_0,\mathbf{x}^N_0 \bigr)$ and take 
$$P_N = \mathcal{L} \bigl( (X^{i,N}, \delta_{\alpha_t^{i,N}}\otimes dt )_{1\leq i \leq N} \bigr) $$
where $\mathcal{L}$ denotes the law under $\mathbb{P}$. This choice of $\mathbb{P}_N$ is admissible for the weak formulation and we have 
$$\mathbb{E}^{P_N} \left [ \int_{t_0}^T \left(  \int_{\R^d}  \frac{1}{N } \sum_{i=1}^N L \bigl(X_t^{i,N}, a \bigr) d\Lambda^{i,N}_t(a) + \mathcal{F} \bigl(\hat{\mu}^N_{t} \bigr) \right)dt  + \mathcal{G} \bigl( \hat{\mu}_T^{N} \bigr) \right] = J^N \bigl(t_0, \mathbf{x}_0;(\alpha_t^{i,N})_{1 \leq i \leq N} \bigr).$$
This implies that
\begin{align*}
\mathcal{V}^N_{w}(t_0, \mathbf{x}_0) & = \inf_{P_N \in \mathcal{R}^N} \mathbb{E}^{P_N} \left [ \int_{t_0}^T \left(  \int_{\R^d}  \frac{1}{N } \sum_{i=1}^N L \bigl(X_t^{i,N}, a \bigr) d\Lambda^{i,N}_t(a) + \mathcal{F} \bigl(\hat{\mu}^N_{t} \bigr) \right)dt  + \mathcal{G} \bigl( \hat{\mu}_T^{N} \bigr) \right] \\
&\leq \inf_{(\alpha_t^{i,N})_{1\leq i \leq N}} J^N \bigl(t_0, \mathbf{x}^N_0;(\alpha_t^{i,N})_{1 \leq i \leq N} \bigr) \\
&= \mathcal{V}^N \bigl(t_0,\mathbf{x}^N_0 \bigr),
\end{align*}
and proves the lemma.
\end{proof}

\subsection{The limit problem}

\label{sec: The mean-field problem 2022}

For some $t_0 \in [0,T]$ and $\mu_0 \in \mathcal{P}_2(\R^d)$ such that $\Psi(\mu_0) <0$, the constrained problem is
\begin{equation}
 \inf_{(\mu, \alpha)} J(t_0,\mu_0; (\mu,\alpha))
 \tag{P}
\label{ConstrainedProblem}
\end{equation}
with $J(t_0,\mu_0;(\mu, \alpha))$ defined by
\begin{equation}
J(t_0,\mu_0;(\mu,\alpha)) := \int_{t_0}^T \int_{\R^d} L \bigl(x,\alpha(t,x) \bigr)d\mu(t)(x) dt + \int_{t_0}^T \mathcal{F}(\mu(t))dt + \mathcal{G}(\mu(T)) 
\label{CostJ11Mai2023}
\end{equation}
and the infimum is taken over couples $(\mu, \alpha) \in \mathcal{C}([t_0, T], \mathcal{P}_2(\R^d)) \times L^2_{dt \otimes \mu(t)} ([t_0,T] \times \R^d, \R^d)$ satisfying in the sense of distributions the Fokker-Planck equation 
\begin{equation}
\left \{
\begin{array}{ll}
 \partial_t \mu + \mdiv(\alpha(t,x) \mu) + \mdiv(b(x,\mu(t))\mu) - \Delta \mu = 0 & \mbox{ in }(t_0,T) \times \R^d \\
\mu(t_0)=\mu_0,
\end{array}
\right.
\label{FPwithdriftb10janv2023}
\end{equation}
under the constraint that $\Psi(\mu(t)) \leq 0$ for all $t \in [t_0,T]$.

The reader may notice that, under appropriate conditions on $\alpha$, the Fokker-Planck equation \eqref{FPwithdriftb10janv2023} describes the law of the solution to the McKean-Vlasov stochastic differential equation
$$dX_t = \alpha(t,X_t)dt + b \bigl(X_t, \mathcal{L}(X_t) \bigr)dt + \sqrt{2}dB_t, \quad X_{t_0} \sim \mu_0, \quad t\in [t_0,T].$$
In this case, the cost \eqref{CostJ11Mai2023} can be rewritten as
$$J(t_0,\mu_0;(\mu,\alpha))= \E \left [ \int_{t_0}^T L \bigl(X_t, \alpha(t,X_t) \bigr) dt + \int_{t_0}^T \mathcal{F} \bigl(\mathcal{L}(X_t) \bigr)dt + \mathcal{G} \bigl(\mathcal{L}(X_T) \bigr) \right], $$
and we could have formulated Problem \eqref{ConstrainedProblem} in terms of optimal control of SDEs of McKean-Vlasov type. As proved in Lemma \ref{lem:weakvsstronglimitproblem}, the resulting value function would be equal to $\mathcal{U}$. We feel that the formulation in terms of pdes is more convenient to address problems with constraints in law.

Of course, Problem \eqref{ConstrainedProblem} is equivalent to the following Problem
\begin{equation}
 \inf_{(\mu,\beta)} J'(t_0,\mu_0; (\mu,\beta))
 \tag{P'}
\label{ConstrainedProblemFormulationbeta}
\end{equation}
with $J' \bigl(t_0,\mu_0;(\mu,\beta) \bigr)$ defined by $ J'\bigl(t_0,\mu_0;(\mu,\beta) \bigr) := J \Bigl(t_0,\mu_0,(\beta - b\bigl(x,\mu(t)),\mu \bigr) \Bigr)$ and the infimum  taken over couples $(\mu, \beta) \in \mathcal{C}([t_0, T], \mathcal{P}_2(\R^d)) \times L^2_{dt \otimes \mu(t)} ([t_0,T] \times \R^d, \R^d)$ satisfying in the sense of distributions the (linear) Fokker-Planck equation 
\begin{equation}
\left \{
\begin{array}{ll}
 \partial_t \mu + \mdiv(\beta(t,x) \mu) - \Delta \mu = 0 & \mbox{ in }(t_0,T) \times \R^d \\
\mu(t_0)=\mu_0,
\end{array}
\right.
\end{equation}
under the constraint that $\Psi(\mu(t)) \leq 0$ for all $t \in [t_0,T]$.
In particular, $\mathcal{U}(t_0,\mu_0)$ can be defined as the infimum in both problems indifferently. Moreover, $(\mu,\alpha)$ is an optimal solution to Problem \eqref{ConstrainedProblem} if and only if $ \bigl(\mu,\alpha + b(x,\mu(t)) \bigr)$ is a solution to Problem \eqref{ConstrainedProblemFormulationbeta}. The advantage of looking at Problem \eqref{ConstrainedProblemFormulationbeta} is that the Fokker-Planck equation reads as a convex constraint in $(\mu,\beta \mu)$. \\

\textbf{Martingale formulation.}  We introduce the controlled martingale formulation. We let $(X,\Lambda)$ be the identity processes over $(\mathcal{C}^d \times \mathcal{V})$ and we look for probabilities $m$ over $\mathcal{C}^d \times \mathcal{V}$ such that $X_0$ is distributed according to $\mu_0$ under $m$, 
$$ \varphi(X_t) - \int_{t_0}^t \int_{\R^d} \mathcal{L}\varphi(X_s \# m, X_s,a) d\Lambda_s(a) ds $$
is a martingale under $m$ for all smooth compactly supported $\varphi : \R^d \rightarrow \R$,
with $\mathcal{L} \varphi (\mu,x,a) = D \varphi (x) . a + D\varphi(x).b(x,\mu)+ \Delta \varphi(x) $. The measure $m$ is also assumed to satisfy the constraint
$$ \Psi(X_t \# m) \leq 0, \quad \forall t \in [t_0,T]. $$
We denote by $\mathcal{R}$ the set of such measures and we look for $m \in \mathcal{R}$ which minimizes the cost function
$$\Gamma(m) := \mathbb{E}^m \left[ \int_{t_0}^T \int_{\R^d} L(X_t, a ) d\Lambda_t(a)dt \right]+ \int_{t_0}^T \mathcal{F}(X_t \#m) dt+ \mathcal{G}(X_T \# m). $$
We denote by $\mathcal{U}_w(t_0,\mu_0)$ the resulting value function.

\begin{lemma}
For all $(t_0,\mu_0) \in [t_0,T] \times \mathcal{P}_2(\R^d)$ such that $\Psi(\mu_0)<0$ it holds
$$\mathcal{U}(t_0,\mu_0) = \mathcal{U}_w(t_0,\mu_0).$$
\label{lem:weakvsstronglimitproblem}
\end{lemma}
\begin{proof}
If we denote by $\mathcal{U}_s(t_0,\mu_0)$ the value of the mean-field problem if defined in terms of controlled stochastic differential equations of McKean-Vlasov type over a fixed probability space, we have, following Theorem 2.4. in \cite{Lacker2017}, $\mathcal{U}_w(t_0,\mu_0) =  \mathcal{U}_s(t_0,\mu_0)$ for all $(t_0,\mu_0) \in [0,T] \times \mathcal{P}_2(\R^d)$ such that $\Psi(\mu_0) \leq 0$. Notably, the presence of the mean-field constraint does not change the argument leading to  the aforementioned result. On the other hand, the optimality conditions of Proposition \eqref{TheoremOptimalityConditions} provide a Lipschitz feedback control for $\mathcal{U}(t_0,m_0)$ when $\Psi(m_0)<0$. We can use this control in the strong formulation to infer that $\mathcal{U}_s(t_0,\mu_0) = \mathcal{U}(t_0,\mu_0)$. This leads to $\mathcal{U}(t_0,\mu_0) = \mathcal{U}_w(t_0,\mu_0)$.
\end{proof}

\subsection{Main result}

Our main result can be stated as follows.

\begin{theorem}
Let Assumption \eqref{AssumptionMeanFieldLimitSept2022} hold. Take $(t_0,\mu_0) \in [0,T] \times \mathcal{P}_2(\R^d)$ such that  $\Psi(\mu_0) < 0$. Then
$$ \lim_{N \rightarrow +\infty} \mathcal{V}^N \bigl(t_0,x_0^{1,N},\dots,x_0^{N,N} \bigl) = \mathcal{U}(t_0,\mu_0), $$
whenever $\frac{1}{N}\sum_{i=1}^N \delta_{x_0^{i,N}} \rightarrow \mu_0$ in $\mathcal{P}_2(\R^d)$ as $N \rightarrow +\infty$.
\label{MainTheoremMFconstraint2022}
Moreover if $P_N$ a sequence of $\epsilon_N$-optimal solutions to the weak $N$-particles problem, for some sequence $\epsilon_N \rightarrow 0$, then the sequence $\hat{\nu}^N \# P_N$ is relatively compact in $\mathcal{P}_p(\mathcal{P}_p(C^d \times \mathcal{V} ))$ for every $p \in (1,2)$. Every limit point is supported on the set of solutions to the mean-field problem in its controlled martingale formulation.
\end{theorem}

A direct consequence of our theorem is that we have a stronger convergence result if the limit problem has a unique solution. In particular, arguing as in Proposition \ref{PropconcentrationAppendix} we can infer that 
$$\lim_{ N\rightarrow +\infty} \E^{P_N} \left[ \sup_{t\in [t_0,T]} d_1 \bigl( \hat{\mu}_t^N , X_t \#m \bigr) \right] =0,$$
if $m$ is the unique solution to the mean field problem in the martingale formulation. Moreover, in this case any solution $(\alpha,\mu)$ of Problem \eqref{ConstrainedProblem} must satisfy $\mu(t) = X_t \#m,$ $\forall t\in [t_0,T]$ and therefore we have
$$\lim_{ N\rightarrow +\infty} \E \left[ \sup_{t\in [t_0,T]} d_1 \bigl( \hat{\mu}_t^N , \mu(t) \bigr) \right] =0.$$

The proof of Theorem \ref{MainTheoremMFconstraint2022} is given in Subsection \ref{sec: from as to mf constraint 11Mai}.

\section{Properties of the mean-field problem}

\label{sec: Properties of the mean-field problem 11 mai 2023}

\subsection{Optimality conditions and regularity of optimal controls}

When there is no constraint, optimal controls for Problem \eqref{ConstrainedProblem} can be characterized by a coupled forward-backward pde system which was first investigated for its connection with mean-field games, see \cite{Briani2018, Lasry2007} for the seminal work of Lasry and Lions and for a derivation of the optimality conditions. In the presence of state constraint, the optimality conditions involve an additional Lagrange multiplier $\nu \in \mathcal{M}^+([0,T])$ which is only active when the optimal trajectory $(\mu(t))_{t\in [t_0,T]}$ touches the constraint. The next result is a small extension of the main result of \cite{Daudin2021}.

\begin{proposition}
Take $(t_0,\mu_0) \in [0,T] \times \mathcal{P}_2(\R^d)$ such that $\Psi(\mu_0)<0$. Under Assumption \eqref{AssumptionMeanFieldLimitSept2022},  Problem \eqref{ConstrainedProblem} admits at least one solution and, for any solution $(\mu,\alpha)$ there exist $u \in L^{\infty}([t_0,T] , E_n ) $ and $\nu \in \mathcal{M}^+([t_0,T])$ such that $\alpha = -\partial_pH(x, Du) $ and $(u,\mu,\nu)$ satisfies the system of optimality conditions
\begin{equation}
\left\{
\begin{array}{lr}
\displaystyle -\partial_t u(t,x) + H(x,Du(t,x)) - Du(t,x).b(x,\mu(t)) - \Delta u(t,x) \\
\displaystyle \hspace{10pt} = \nu(t) \frac{\delta \Psi}{\delta m}(\mu(t),x) + \frac{\delta \mathcal{F}}{\delta m}(\mu(t),x) + \int_{\R^d} Du(t,y). \frac{\delta b}{\delta m}(y, \mu(t),x)d\mu(t)(y) & \mbox{in   } (t_0,T)\times \R^d, \\
\displaystyle \partial_t \mu - \mdiv( \partial_pH(x, Du(t,x)) \mu ) +\mdiv(b(x,\mu(t))\mu)- \Delta \mu = 0  & \mbox{in   } (t_0,T)\times \R^d,  \\
\displaystyle \mu(t_0) = \mu_0, \hspace{30pt} \displaystyle u(T,x) = \frac{\delta \mathcal{G}}{\delta m}(\mu(T),x)  & \mbox{       in      } \R^d,
\end{array}
\right.
\label{OptimalityConditionMainTheorem2021}
\end{equation}
where the Fokker-Planck equation is understood in the sense of distributions and $u$ solves the HJB equation in the sense of Definition \ref{definitionHJBsol2023} below.  
The Lagrange multiplier $\nu$ satisfies the exclusion condition
\begin{equation}
\Psi(\mu(t)) = 0, \hspace{50pt} \mbox{$\nu$-almost-everywhere in }[t_0,T] .
\label{Exclusionnu}
\end{equation}
The optimal control $\alpha$ belongs to $BV_{loc}([t_0,T] \times \R^d , \R^d) \bigcap L^{\infty}([t_0,T], \mathcal{C}_b^{n-1}(\R^d,\R^d))$.
Finally the value of the optimal control problem is given by
$$\mathcal{U}(t_0,\mu_0) = \int_{\R^d} u(t_0,x)d\mu_0(x) + \int_{t_0}^T \mathcal{F}(\mu(t))dt + \mathcal{G}(\mu(T)).$$
\label{TheoremOptimalityConditions}
\end{proposition}

\begin{remark}
Under additional regularity and qualification conditions on $\Psi$ ---Assumptions \eqref{APsiC2} and \eqref{TransCondPsiMF2022} in Section \ref{sec: Application to Large Deviations Sept2022}---, it is shown in \cite{Daudin2021} that $\nu = \nu_1 + \eta \delta_T$ for some $\nu_1 \in L^{\infty}([0,T])$ and some $\eta \geq 0$. As a consequence $u$ belongs to $\mathcal{C}([0,T],E_n)$ and optimal controls are Lipschitz continuous in time and space. See Theorem 2.2 in \cite{Daudin2021}.
\end{remark}

In the above theorem we use the following notion of solution for the HJB equation, using Duhamel's representation formula.
\begin{definition}
Let $\psi_1,\varphi_1 \in \mathcal{C}^0([t_0,T],E_n)$, $\psi_2 \in E_{n+\alpha}$ for some $n \geq 2$ and $\nu \in \mathcal{M}^+([0,T])$. We say that $u \in L^{\infty}([t_0,T], E_n)$ is a solution to 
\begin{equation}
\left \{
\begin{array}{ll}
-\partial_t u + H(x,Du) - \Delta u = \psi_1 \nu + \varphi_1 & \mbox{ in }[t_0,T] \times \R^d,\\
u(T,x) = \psi_2 & \mbox{ in } \R^d,
\end{array}
\right.
\end{equation}
if, for almost all $t \in [t_0,T]$, 
\begin{align}
\notag u(t,x) &= P_{T-t}\psi_2(x) +\int_{t_0}^T \mathds{1}_{(t,T]}P_{s-t}\psi_1(s,.)(x)d\nu(s) + \int_t^T P_{s-t} \varphi_1(s,.)(x)ds \\
&- \int_{t}^T P_{s-t}\left[H(.,Du(s,.))\right] (x)ds, \hspace{30pt} \forall x\in \R^d, 
\label{RepresentationformulaHJB10Mai}
\end{align}
where $(P_t)_{t \geq 0}$ is the heat semi-group (associated to $-\Delta$).
\label{definitionHJBsol2023}
\end{definition}
This formulation is convenient to handle solutions which are not necessarily continuous in time (the optimal control can jump when the optimal trajectory touches the constraint) but are, for each time, regular in the space variable.

The proof of Theorem \ref{TheoremOptimalityConditions} is given in \cite{Daudin2021} when $b=0$ by a penalization procedure. We give here a more direct proof covering the case $b \neq 0$ using a min/max argument similar to \cite{Daudin2022} Section 3.2. See  Appendix \ref{subsec:OptimalityConditions}. 

\subsection{Stability with respect to the constraint}

\label{sec: Stability section}

For $\delta \geq 0$ small, we define $\mathcal{U}^{\delta}(t_0,\mu_0)$ to be the value of the same problem associated to the constraint $\Psi(\mu(t)) \leq -\delta$ for all $t \in [t_0,T]$. In particular it holds that $\mathcal{U}^{\delta_1}(t_0,\mu_0) \geq \mathcal{U}^{\delta_2}(t_0,\mu_0)$ whenever $\delta_1 \geq \delta_2 \geq 0$. Using the convexity of the constraint, we can prove the following stability result.

\begin{lemma} 
Let Assumption \eqref{AssumptionMeanFieldLimitSept2022}  hold and assume as well that $\Psi(\mu_0) <0$. Then it holds
$$ \lim_{\delta \rightarrow 0} \mathcal{U}^{\delta}(t_0,\mu_0) = \mathcal{U}(t_0,\mu_0).$$
\label{StabilityLemma7octobre}
\end{lemma}
\begin{proof}
Assume on the contrary that 
\begin{equation} 
\lim_{\delta \rightarrow 0} \mathcal{U}^{\delta}(t_0,\mu_0) = \inf_{\delta >0} \mathcal{U}^{\delta}(t_0,\mu_0) = \mathcal{U}(t_0,\mu_0)+ \gamma
\label{towardacontradiction5avril2023}
\end{equation}
for some $\gamma >0$. For every $\delta >0$ we denote by $(\mu^{\delta}, \beta^{\delta})$ an optimal solution for $\mathcal{U}^{\delta}(t_0,\mu_0)$ and by $(\tilde{\beta},\tilde{m})$ and optimal solution for $\mathcal{U}(t_0,\mu_0)$ where both Problems are understood in term of Problem \eqref{ConstrainedProblemFormulationbeta}. For $\lambda \in (0,1)$ we let $\mu^{\delta,\lambda} := (1-\lambda) \tilde{\mu} + \lambda \mu^{\delta}$ and, noticing that $(1-\lambda) \tilde{\beta} \tilde{\mu} + \lambda \beta^{\delta} \mu^{\delta}$ is absolutely continuous with respect to $\mu^{\delta,\lambda}$, we  define a control $\beta^{\delta, \lambda}$ such that $(1-\lambda) \tilde{\beta} \tilde{\mu} + \lambda \beta^{\delta} \mu^{\delta} = \beta^{\delta, \lambda} \mu^{\delta, \lambda}$, so that $(\mu^{\delta,\lambda})_{t\in [t_0,T]}$ satisfies
$$ \partial_t \mu^{\delta, \lambda} + \mdiv(\beta^{\delta \lambda} \mu^{\delta, \lambda}) -\Delta \mu^{\delta, \lambda} = 0 \quad \mbox{ in } (t_0,T) \times \R^d, \hspace{30pt} \mu^{\delta,\lambda}(t_0) = \mu_0. $$
By convexity of $\Psi$, we have, for all $t\in [t_0,T]$,
$$ \Psi( \mu^{\delta, \lambda}(t)) \leq \lambda \Psi(\mu^{\delta}(t)) \leq -\delta \lambda. $$
As a consequence, $(\beta^{\delta,\lambda}, \mu^{\delta,\lambda})$ is admissible for $\mathcal{U}^{\delta\lambda}(t_0,\mu_0)$ and 
\begin{equation}
J'(t_0,\mu_0,(\beta^{\delta,\lambda}, \mu^{\delta,\lambda}) \geq \mathcal{U}^{\delta \lambda}(t_0,\mu_0).
\label{togetherwithleadstocontradiction}
\end{equation}
However, for all $\lambda \in (0,1)$ it holds, by convexity of $\R^d \times \R^+ \ni (a,b) \mapsto L(x,\frac{a}{b})b$ (set to be $+\infty$ if $b=0$),
\begin{align*}
J'(t_0,\mu_0,(\beta^{\delta,\lambda}, \mu^{\delta,\lambda})) &\leq (1-\lambda) \int_{t_0}^T \int_{\R^d} L(x, \tilde{\beta}(t,x) - b(x, \mu^{\delta,\lambda}(t)) d\tilde{\mu}(t)(x)dt \\
&+ \lambda \int_{t_0}^T \int_{\R^d} L(x, \beta^{\delta}(t,x) - b(x, \mu^{\delta,\lambda}(t))d\mu^{\delta}(t))dt \\
&+ \int_{t_0}^T \mathcal{F}(\mu^{\delta,\lambda}(t))dt + \mathcal{G}(\mu^{\delta,\lambda}(T)).
\end{align*}
This shows in particular that $\limsup_{\lambda \rightarrow 0^+} J'(t_0,\mu_0,(\beta^{\delta,\lambda}, \mu^{\delta,\lambda})) \leq J'(t_0,\mu_0,(\tilde{\beta},\tilde{\mu}))=\mathcal{U}(t_0,\mu_0) $. Together with \eqref{togetherwithleadstocontradiction} this contradicts \eqref{towardacontradiction5avril2023}.
\end{proof}

\begin{remark}
We easily conclude from the above lemma that the value of the limit problem does not change if we replace the closed constraint $\{ \Psi \leq 0 \}$ by the open one $\{ \Psi <0 \}$.
\end{remark}

\section{Mean field limit}

\label{sec: Mean field limit 2022}

The main result of this section is to prove Theorem \eqref{MainTheoremMFconstraint2022}, that is the convergence of $\mathcal{V}^N(t_0,\mathbf{x}_0)$ to  $\mathcal{U}(t_0,\mu_0)$ as $\displaystyle \frac{1}{N}\sum_{i=1}^N \delta_{x_0^{i,N}} \rightarrow \mu_0$ when $N \rightarrow +\infty$.

\subsection{From mean-field to almost-sure constraint}

In this section we prove the first inequality in Theorem \eqref{MainTheoremMFconstraint2022}.

\begin{proposition}
Let Assumption \eqref{AssumptionMeanFieldLimitSept2022} hold. Assume further that $\mu_0$ satisfies $\Psi(\mu_0)<0$. Then it holds that
$$ \limsup_{ N \rightarrow +\infty} \mathcal{V}^N \bigl(t_0, \mathbf{x}_0^N \bigr) \leq \mathcal{U}(t_0,\mu_0), $$
whenever $\lim_{N \rightarrow +\infty} \frac{1}{N}\sum_{i=1}^N \delta_{x_0^{i,N}} = \mu_0$ in $\mathcal{P}_2(\R^d)$.
\label{TheoremFirstInequalitySept2022}
\end{proposition}

\begin{proof}

To prove Proposition \eqref{TheoremFirstInequalitySept2022} we proceed as follows. First we fix $\delta \in (0,-\Psi(\mu_0)/2)$ small and we take $\alpha : [t_0,T] \times \R^d \rightarrow \R^d$ to be an optimal control for $\mathcal{U}^{\delta}(t_0,\mu_0)$. Using Theorem \ref{TheoremOptimalityConditions}, we know that $\alpha$ is bounded and Lipschitz continuous in space. We let $\mu$ be the corresponding trajectory, solution to
\begin{equation*}
\left \{
\begin{array}{ll}
\partial_t \mu + \mdiv(\alpha \mu)+ \mdiv(b\mu)- \Delta \mu = 0 & \mbox{ in }(t_0,T) \times \R^d, \\
\mu(t_0)=\mu_0.
\end{array}
\right.
\end{equation*}
In particular, $\Psi(\mu(t)) \leq -\delta$ for all $t\in [t_0,T]$. For a set of initial positions $\mathbf{x}_0^N = \bigl(x_0^{1,N},\dots,x_N^{N,N} \bigr) \in \bigl(\R^d \bigr)^N$ such that $\frac{1}{N}\sum_{i=1}^N \delta_{x_0^{i,N}} \rightarrow \mu_0$ in $\mathcal{P}_2(\R^d)$,  we let $ \bigl(X_t^{1,N}, \dots , X_t^{N,N} \bigr)_{t_0 \leq t \leq T}$ be the solution to 
\begin{equation*}
\left \{
\begin{array}{ll}
\displaystyle X_t^{i,N} = x_0^{i,N} + \int_{t_0}^{t \wedge \tau_N} \alpha \bigl(s,X_s^{i,N} \bigr) dt + \int_{t \wedge \tau_N }^t \beta_s^{i,N} ds + \int_{t_0}^t b \bigl(X_s^{i,N}, \hat{\mu}_s^{N} \bigr)ds+ \sqrt{2} \int_{t_0}^t dB_s^{i,N}, \quad 1\leq i \leq N  \\
\displaystyle \hat{\mu}_s^N = \frac{1}{N} \sum_{i=1}^N \delta_{X_s^{i,N}}
 \end{array}
\right.
\end{equation*}
where 
$$\tau_N := \inf \bigl \{t \geq 0, \Psi \bigl( \hat{\mu}_t^{N} \bigr ) \geq -\frac{\delta }{2} \bigr \},$$ 
with the convention $\inf \{ \emptyset \} = +\infty$, and $\beta_t^{i,N}$ is the feedback control, defined for all $t \geq \tau_N \wedge T$ by 
$$\beta_t^{i,N} =  \frac{4 \bigl(X_t^{i,N} - X_{\tau_N\wedge T}^{i,N} \bigr)}{\sum_{i=1}^N | X_t^{i,N} - X_{\tau_N\wedge T}^{i,N} |^2 -r^2 N} - 2\frac{d}{r^2} \bigl(X_t^{i,N}-X_{\tau_N\wedge T}^{i,N} \bigr) -b \bigl(X_t^{i,N}, \hat{\mu}_t^{N} \bigr), $$
with $r = \frac{\delta}{4C_{\Psi}}$ and $C_{\Psi}$ a Lipschitz constant for $\Psi$ with respect to $d_1$. We also assume that $N$ is large enough so that $\Psi \Bigl( \frac{1}{N}\sum_{i=1}^N \delta_{x_0^{i,N}} \Bigr) < -\delta $. We will need the following key lemma, which justifies how we chose $\bigl (\beta^{i,N})_{1\leq i \leq N}$.

\begin{lemma}
$\mathbb{P}$-almost-surely, it holds that,
$$ \frac{1}{N} \sum_{i=1}^N \left| X_t^{i,N} - X_{\tau_N\wedge T}^{i,N} \right|^2 \leq r^2, \hspace{30pt} \forall t \geq \tau_N\wedge T. $$
Moreover, the following estimate holds  
\begin{align*}
\E \left[ \int_{\tau_N \wedge T}^{T} \frac{1}{N} \sum_{i=1}^N \left| \beta_t^{i,N} \right|^2 dt \right] & \leq \frac{32 d}{r^2N} \E \left[e^{T- T \wedge \tau_N} \right] + \frac{16d^2}{r^2} \E \left[T- T\wedge \tau_N \right] \\
&+ 2 \|b\|^2_{\infty} \E \left[T- T\wedge \tau_N \right]  .
\end{align*}
\label{FreezeLemma18Octobre2022}
\end{lemma}

We continue with the ongoing proof. We have taken $r$ and $\beta_N$ is such a way that $\mathbb{P}$-almost-surely $\displaystyle \Psi( \hat{ \mu}_t^{N} ) \leq \frac{-\delta}{4}$ for all  $t \in [t_0,T]$. Indeed, by definition of $\tau_N$, $\mathbb{P}$-almost-surely, $ \Psi( \hat{\mu}_t^{N} ) \leq -\frac{\delta}{2}$ for all $ t \leq \tau_N,$  and, $\mathbb{P}$-almost-surely, by definition of $r$ and Lemma \ref{FreezeLemma18Octobre2022}, it holds, whenever $t \geq \tau_N$
\begin{align*}
\left | \Psi(\hat{\mu}_t^{N}) - \Psi ( \hat{\mu}_{\tau_N}^{N}) \right| &\leq C_{\Psi} d_1 \bigl( \hat{\mu}_t^{N}, \hat{\mu}_{\tau_N}^{N} \bigr) \\
&  \leq C_{\Psi} d_2 \bigl( \hat{\mu}_t^{N}, \hat{\mu}_{\tau_N}^{N} \bigr) \\
&\leq \frac{C_{\Psi}}{\sqrt{N}} \left( \sum_{i=1}^N |X_t^{i,N} - X_{\tau_N}^{i,N} |^2 \right)^{1/2} \\
&\leq C_{\Psi} \times r \leq \frac{\delta}{4}
\end{align*}
and, as a consequence, being $\displaystyle \Psi \bigl(\hat{\mu}^{N}_{\tau_N} \bigr) = -\delta/2$, it holds that $\Psi \bigl(\hat{\mu}_t^{N} \bigr) \leq -\delta/4$, for all $t \in [t_0,T]$. Therefore, we have an admissible control for $\mathcal{V}^N \bigl(t_0,\mathbf{x}_0^N \bigr)$. Now, by standard propagation of chaos estimates, see Proposition \eqref{PropconcentrationAppendix} in Apendix \ref{sec: Concentration limit}, it holds that 
$$ \lim_{N \rightarrow +\infty} \E \left[ \sup_{t \in [t_0,T\wedge \tau_N]} d_1 \bigl(\mu(t), \hat{\mu}^{N}_t \bigr) \right]  =0.$$
As a consequence, using that $\Psi \bigl(\mu(t) \bigr) \leq -\delta$, for all $t \in [t_0,T]$ as well as the Lipschitz continuity of $\Psi$ with respect to $d_1$, we get
\begin{align*}
\mathbb{P} \left[ \tau_N < T \right] & = \mathbb{P} \left[ \exists t <T, \Psi \bigl( \hat{\mu}^{N}_t \bigr ) \geq \frac{-\delta}{2} \right] \\
& \leq  \mathbb{P} \left[ \exists t <\tau_N, \Psi \bigl( \hat{\mu}^{N}_t \bigr) \geq \frac{-3\delta}{4} \right] \\
& \leq \mathbb{P} \left[ \sup_{t \in [t_0,T \wedge \tau_N] } d_1 \bigl( \hat{\mu}^{N}_t, \mu(t) \bigr) \geq \frac{\delta}{4C_{\Psi}} \right].
\end{align*}
Using Markov's inequality we conclude that
\begin{align*}
\mathbb{P} \left[ \tau_N < T \right]& \leq \frac{4C_{\Psi}}{\delta} \E \left[\sup_{t \in [t_0,T \wedge \tau_N] } d_1 \bigl( \hat{\mu}^{N}_t, \mu(t) \bigr) \right].
\end{align*}
Being $ \E  \left[ T-T\wedge \tau_N \right] \leq T \mathbb{P} \left[ \tau_N <T \right]$ we conclude that $\displaystyle  \lim_{N \rightarrow +\infty} \E \left[ T- T \wedge \tau_N \right ] =0.$ Now we can use Lemma \eqref{FreezeLemma18Octobre2022} and get $\displaystyle \lim_{N \rightarrow +\infty} \E^{\mathbb{P}^{\gamma_N}} \left[ \int_{ \tau_N \wedge T}^T \frac{1}{N} \sum_{i=1}^N |\beta_t^{i,N} |^2 \right] =0, $ and we easily deduce
$$\displaystyle \lim_{N \rightarrow +\infty}  \E \left [ \sup_{t \in [t_0,T]} d_1 (\mu(t), \hat{\mu}^{N}_t) \right]  =0.$$
As a consequence, $\alpha$ being bounded and Lipschitz continuous in the space variable and $\mathcal{F}$ and $\mathcal{G}$ being Lipschitz continuous with respect to $d_1$, 
\begin{align*}
& \lim_{N \rightarrow +\infty}   \E \Bigl [ \int_{t_0}^{T \wedge \tau_N} \frac{1}{N} \sum_{i=1}^N L \Bigl(X_t^{i,N}, \alpha \bigl(t,X_t^{i,N} \bigr) \Bigr) dt + \int_{t_0}^{T \wedge \tau_N} \mathcal{F} \bigl(\hat{\mu}_t^{N} \bigr) dt  \\
&+  \int_{T \wedge \tau_N}^T \frac{1}{N} \sum_{i=1}^N L \Bigl(X_t^{i,N}, \beta_t^{i,N}) \Bigr) dt + \int_{T \wedge \tau_N}^T \mathcal{F} \bigl(\hat{\mu}_t^{N} \bigr) dt  + \mathcal{G} \bigl(\hat{\mu}_T^{N,x} \bigr) \Bigr ] \\
&= \lim_{N \rightarrow +\infty} \E \Bigl [ \int_{t_0}^{T\wedge \tau_N} \int_{\R^d} L \bigl(x,\alpha(t,x) \bigr)d\hat{\mu}_t^N(x)dt + \int_{t_0}^{T \wedge \tau_N} \mathcal{F} \bigl(\hat{\mu}_t^{N} \bigr) dt  \\
&+ \int_{T \wedge \tau_N}^T \frac{1}{N} \sum_{i=1}^N L \big( X_t^{i,N}, \beta_t^{i,N} \bigr) dt + \int_{T \wedge \tau_N}^T \mathcal{F} \bigl(\hat{\mu}_t^{N} \bigr) dt   + \mathcal{G} \bigl(\hat{\mu}_T^{N,x} \bigr)  \Bigr ] \\ 
&= \int_{t_0}^T \int_{\R^d} L \bigl(x,\alpha(t,x) \bigr) d\mu(t)(x)dt + \int_{t_0}^T \mathcal{F} \bigl(\mu(t) \bigr) dt +\mathcal{G}\bigl(\mu(T) \bigr).
\end{align*}
Finally, being $\alpha$ optimal for $\mathcal{U}^{\delta}(t_0,\mu_0)$ we have that 
$$\limsup_{N \rightarrow + \infty} \mathcal{V}^N \bigl(t_0,\mathbf{x}_0^N \bigr) \leq \mathcal{U}^{\delta}(t_0,\mu_0).$$
Yet, we have proved, in Proposition \eqref{StabilityLemma7octobre} that  $\lim_{\delta \rightarrow 0}  \mathcal{U}^{\delta}(t_0,\mu_0) = \mathcal{U}(t_0,\mu_0)$ and therefore,
$$\limsup_{N \rightarrow + \infty} \mathcal{V}^N \bigl(t_0,\mathbf{x}_0^N \bigr) \leq \mathcal{U}(t_0,\mu_0),$$
which concludes the proof of the proposition.
\end{proof}

It remains to prove Lemma \ref{FreezeLemma18Octobre2022}.

\begin{proof}[Proof of Lemma \ref{FreezeLemma18Octobre2022} ]
For $\eta \geq0$ small, we introduce the stopping time 
$$ \tau^{\eta} := \inf \Bigl \{ t \geq \tau_N\wedge T, \frac{1}{N} \sum_{i=1}^N |X_t^{i,N} - X_{\tau_N \wedge T}^{i,N}|^2 \geq r^2-\eta \Bigr \},$$
with the convention that $\inf \left \{ \emptyset \right \} = +\infty$. 

For $\eta >0$ and  $T'>T$, we write $\textbf{B}_t = ^t\Bigl(B_t^{1,N}, \dots, B_t^{N,N} \Bigr)$ and $\textbf{Y}_t = ^t \Bigl(X_t^{i,N}, \dots, X_t^{N,N} \Bigr)$ and apply Itô's lemma to get 
\begin{align*}
-e^{-\tau^{\eta}\wedge T'} &\log \Bigl(r^2 - \frac{|\textbf{Y}_{\tau^{\eta}\wedge T'}-\textbf{Y}_{\tau_N\wedge T}|^2}{N} \Bigr) = -e^{-\tau_N\wedge T} \log(r^2) + \int_{\tau_N\wedge T}^{\tau^{\eta}\wedge T'} e^{-t} \log \Bigl( r^2 - \frac{|\textbf{Y}_t-\textbf{Y}_{\tau_N\wedge T}|^2}{N} \Bigr) dt \\
&+ \int_{\tau_N\wedge T}^{\tau^{\eta}\wedge T'} e^{-t} \left[ \frac{4 (\textbf{Y}_t - \textbf{Y}_{\tau_N\wedge T})}{|\textbf{Y}_t - \textbf{Y}_{\tau_N\wedge T}|^2-r^2N}.\frac{2(\textbf{Y}_t - \textbf{Y}_{\tau_N\wedge T})}{Nr^2-|\textbf{Y}_t - \textbf{Y}_{\tau_N\wedge T}|^2} - \frac{2d}{r^2} \frac{2|\textbf{Y}_t - \textbf{Y}_{\tau_N\wedge T}|^2}{Nr^2-|\textbf{Y}_t - \textbf{Y}_{\tau_N\wedge T}|^2} \right] dt \\
&+ \int_{\tau_N\wedge T}^{\tau^{\eta}\wedge T'} e^{-t} \left[ \frac{2dN}{Nr^2 -|\textbf{Y}_t - \textbf{Y}_{\tau_N\wedge T}|^2} + \frac{4|\textbf{Y}_t - \textbf{Y}_{\tau_N\wedge T}|^2}{(Nr^2 - |\textbf{Y}_t - \textbf{Y}_{\tau_N\wedge T}|^2)^2} \right] dt \\
& +\sqrt{2} \int_{\tau_N\wedge T}^{\tau^{\eta}\wedge T'} \frac{2 (\textbf{Y}_t - \textbf{Y}_{\tau_N\wedge T})}{Nr^2-|\textbf{Y}_t - \textbf{Y}_{\tau_N\wedge T}|^2}.d\textbf{B}_t \\
&\leq -4 \int_{\tau_N\wedge T}^{\tau^{\eta}\wedge T'} e^{-t} \frac{|\textbf{Y}_t - \textbf{Y}_{\tau_N\wedge T}|^2}{(|\textbf{Y}_t - \textbf{Y}_{\tau_N\wedge T}|^2-r^2N)^2}dt + \int_{\tau_N\wedge T}^{\tau^{\eta}\wedge T'} e^{-t} \left[ -\frac{4d}{r^2} \frac{|\textbf{Y}_t - \textbf{Y}_{\tau_N\wedge T}|^2-dN}{Nr^2- |\textbf{Y}_t - \textbf{Y}_{\tau_N\wedge T}|^2} \right] dt \\
&+ \sqrt{2}\int_{\tau_N\wedge T}^{\tau^{\eta}\wedge T'} \frac{2 (\textbf{Y}_t - \textbf{Y}_{\tau_N\wedge T})}{Nr^2-|\textbf{Y}_t - \textbf{Y}_{\tau_N\wedge T}|^2}.d\textbf{B}_t
\end{align*}
However, an elementary analysis reveals that
$$  -\frac{4d}{r^2} \frac{|\textbf{Y}_t - \textbf{Y}_{\tau_N\wedge T}|^2-dN}{Nr^2- |\textbf{Y}_t - \textbf{Y}_{\tau_N\wedge T}|^2} \leq \frac{2d}{r^2} $$
whenever $0 \leq |\textbf{Y}_t - \textbf{Y}_{\tau_N\wedge T}|^2 < Nr^2$.
Therefore, we get, multiplying by $e^{\tau_N\wedge T}$ and taking expectations,
\begin{align*}
 -  \E & \left[ e^{\tau_N\wedge T-\tau^{\eta}\wedge T'} \log \Bigl(r^2 - \frac{|\textbf{Y}_{\tau^{\eta}\wedge T'}-\textbf{Y}_{\tau_N\wedge T}|^2}{N} \Bigr) \right ] \\
 & + 4 \E\left[ \int_{\tau_N\wedge T}^{\tau^{\eta}\wedge T'} e^{\tau_N\wedge T-t} \frac{|\textbf{Y}_t - \textbf{Y}_{\tau_N\wedge T}|^2}{(|\textbf{Y}_t - \textbf{Y}_{\tau_N\wedge T}|^2-r^2N)^2} dt \right] \leq \frac{2d}{r^2}.
 \end{align*}
Letting $T' \rightarrow +\infty$, using the definition of $\tau^{\eta}$ and Lebesgue dominated convergence theorem leads to 
 \begin{equation}
- \log(\eta) \E \left[ e^{\tau_n\wedge T - \tau^{\eta}} \mathds{1}_{\left \{ \tau^{\eta} < +\infty \right \} } \right] + 4  \E \left[ \int_{\tau_N\wedge T}^{\tau^{\eta}} e^{\tau_N\wedge T-t} \frac{|\textbf{Y}_t - \textbf{Y}_{\tau_N\wedge T}|^2}{(|\textbf{Y}_t - \textbf{Y}_{\tau_N\wedge T}|^2-r^2N)^2} dt \right] \leq \frac{2d}{r^2}.
\label{JEANAYMAROCT2022}
\end{equation}
Notice that both terms in the left-hand side of \eqref{JEANAYMAROCT2022} are non-negative for $\eta \leq 1$. Letting $\eta \rightarrow 0$, we get, on the one hand that  $\tau^0 = +\infty$, $\mathbb{P}$-almost surely and, on the other hand, we obtain
$$4 \E\left[ \int_{\tau_N \wedge T}^{+\infty} e^{\tau_N\wedge T-t} \frac{|\textbf{Y}_t - \textbf{Y}_{\tau_N\wedge T}|^2}{(|\textbf{Y}_t - \textbf{Y}_{\tau_N\wedge T}|^2-r^2N)^2} dt \right] \leq \frac{2d}{r^2}. $$
It follows that, 
\begin{align*}
\E &\left[\int_{\tau_N \wedge T}^T \frac{1}{N} \sum_{i=1}^N |\beta_t^{i,N} |^2 dt \right] \\
& \leq 2\E \left[ \int_{\tau_N \wedge T}^T  \frac{1}{N}  \left| \frac{4(\textbf{Y}_t - \textbf{Y}_{\tau_N})}{|\textbf{Y}_t - \textbf{Y}_{\tau_N}|^2-r^2N} -2\frac{d}{r^2} (\textbf{Y}_t - \textbf{Y}_{\tau_N}) \right|^2 dt \right] +2 \|b\|^2_{\infty} \E \left[ T- T\wedge \tau_N \right] \\
&\leq 4 \E \left[  \int_{\tau_N \wedge T}^T  \frac{1}{N}  \left| \frac{4(\textbf{Y}_t - \textbf{Y}_{\tau_N})}{|\textbf{Y}_t - \textbf{Y}_{\tau_N}|^2-r^2N}  \right|^2 dt \right]  + 4\E^{\mathbb{P}^{\gamma_N}}  \left[  \int_{\tau_N \wedge T}^T   \frac{1}{N}  \left| 2\frac{d}{r^2} (\textbf{Y}_t - \textbf{Y}_{\tau_N}) \right|^2 dt \right]  \\
&+2 \|b\|^2_{\infty} \E \left[ T- T\wedge \tau_N \right] \\
&\leq \frac{32d^2}{r^2N} \E \left[ e^{T- T\wedge \tau_N} \right] + \frac{16d^2}{r^4} \E \left[  \int_{\tau_N \wedge T}^T \frac{1}{N}  |\textbf{Y}_t - \textbf{Y}_{\tau_N} |^2 dt \right] +2 \|b\|^2_{\infty} \E \left[ T- T\wedge \tau_N \right] \\
&\leq \frac{16d^2}{r^2N} \E \left[ e^{T - T\wedge \tau_N} \right] + \frac{8d^2}{r^2} \E \left[T -T\wedge \tau_N \right] +2 \|b\|^2_{\infty} \E \left[ T- T\wedge \tau_N \right],
\end{align*}
where we used, for the last inequality, the fact that, $\mathbb{P}$-almost-surely, for all $t\geq \tau_N$, 
$$ \left|\textbf{Y}_t - \textbf{Y}_{\tau_N} \right|^2 \leq N r^2. $$
This concludes the proof of the lemma.
\end{proof}

\subsection{From almost-sure to mean-field constraint}

\label{sec: from as to mf constraint 11Mai}

To prove the second inequality we rely on compactness methods developed, in the context of Large Deviations by Budhiraja, Dupuis and Fischer \cite{Budhiraja2012} and, in the context of mean-field control, by Lacker \cite{Lacker2017} and Djete, Possamaï and Tan \cite{Djete2020}. 

Recall that we introduced the weak formulation of the $N$-particle problem in Subsection \ref{sec: The N-state problem} and the controlled-martingale formulation of the mean-field problem in Subsection \ref{sec: The mean-field problem 2022}.

Theorem \ref{MainTheoremMFconstraint2022} is equivalent to the following proposition.
\begin{proposition}
Let us fix $\mu_0 \in \mathcal{P}_2(\R^d)$ such that $\Psi(\mu_0) <0$ as well as some $\mathbf{x}^N_0 = (x_0^{1,N},\dots,x_0^{N,N}) \in \Omega_N$ such that $\displaystyle \frac{1}{N}\sum_{i=1}^N \delta_{x_0^{i,N}} \rightarrow \mu_0$ in $\mathcal{P}_2(\R^d)$. Take $P_N$ a sequence of $\epsilon_N$-optimal solutions to the weak $N$-particle problem, for some sequence $\epsilon_N \rightarrow 0$. Then the sequence $\hat{\nu}^N \# P_N$ is relatively compact in $\mathcal{P}_p(\mathcal{P}_p(C^d \times \mathcal{V} ))$ for every $p \in (1,2)$. Every limit point is supported on the set of solutions to the relaxed mean-field problem and it holds that
$$ \mathcal{U}(t_0,\mu_0) = \lim_{ N\rightarrow +\infty} \mathcal{V}^N(t_0,\mathbf{x}_0^N). $$

\label{EasyInequalitySeptember2022}

\end{proposition}

\begin{proof}
We will closely follow the steps of \cite{Lacker2017} and therefore we only highlight the differences due to the constraint. In light of \cite{Lacker2015} Corollary B.2, to prove the pre-compactness of  $\hat{\nu}^N \# P_N$, it suffices to prove that the mean measures $ \displaystyle \frac{1}{N} \sum_{i=1}^N (X^{i,N},\Lambda^{i,N}) \# P_N$ are tight and to prove that
\begin{equation}
\sup_{N} \E^{P_N} \frac{1}{N} \sum_{i=1}^N \left[ \sup_{t\in [t_0,T]} |X_t^{i,N}|^2 + \int_{t_0}^T \int_{\R^d} |a|^2 d\Lambda_t^{i,N}(a)dt \right] <+ \infty.
\label{TightnessalaLacker2022}
\end{equation}
The tightness of the mean measures actually follows from \eqref{TightnessalaLacker2022} thanks to the compactness result of Proposition 3.5 in \cite{Lacker2017}. By standard estimates, it is enough to prove
$$ \sup_{N}  \E^{P_N} \left[ \frac{1}{N}\sum_{i=1}^N |x_0^{i,N}|^2 \right] < +\infty$$
as well as 
$$  \sup_{N} \E^{P_N} \left[\int_{t_0}^T \frac{1}{N}\sum_{i=1}^N \int_{\R^d} |a|^2 d\Lambda_t^{i,N}(a)dt \right] <+ \infty $$
in order to get \eqref{TightnessalaLacker2022}.
The former follows from the convergence in $\mathcal{P}_2(\R^d)$ of $\frac{1}{N} \sum_{i=1}^N \delta_{x_0^{i,N}}$ toward $\mu_0$. The latter follows from the coercivity of $L$, the boundness of $\mathcal{F}$, $\mathcal{G}$ and the fact that we took the $P_N$ as $\epsilon_N$-optimal solutions for the $N$-particle problem whose values are bounded independently from $N$ (as can be deduced from Theorem \eqref{TheoremFirstInequalitySept2022}). Now we take a limit point $P \in \mathcal{P}_p(\mathcal{P}_p(\mathcal{C}^d \times \mathcal{V} ))$ and prove that $P$ is supported on the set of solutions to the mean field relaxed problem. First we have that $\hat{\mu}_0^{N}\# P_N \rightarrow \delta_{\mu_0}$ in $\mathcal{P}_p(\mathcal{P}_p(\R^d))$.  Following \cite{Lacker2017} Proposition 5.2 we have that $P$ is supported on the set of measures solution to the martingale problem. It remains to show that the constraint is satisfied $P$-almost surely at the limit. By continuity of $\Psi$, for all $t\in [t_0,T]$ it holds that 
$$P \left( \left \{ m \in \mathcal{P}_p(\mathcal{C}^d \times \mathcal{V}),  \quad\Psi(X_t \# m) \leq 0) \right \} \right) \geq \limsup_{N \rightarrow +\infty} P_N(\Psi(\hat{\mu}_t^{N}) \leq 0) = 1.$$
Since $P$-almost surely $m$ satisfies the martingale problem, we have that $P$-almost surely $t \rightarrow X_t \# m$ is continuous and therefore we have that
$$P \left( \left \{ m \in \mathcal{P}_p(\mathcal{C}^d \times \mathcal{V}), \quad \Psi(X_t \# m) \leq 0 \quad \forall t\in [t_0,T]) \right \} \right) =1.$$
This implies the $P$ is supported on $\mathcal{R}$, the set of admissible candidate for the mean-field problem, as defined in Subsection \ref{sec: The mean-field problem 2022}. The fact that $P$ is supported on the set of optimal solutions of the mean-field problem follows from the lower semi-continuity of the cost functional as proved in \cite{Lacker2017} Lemma 4.1 and from Proposition \eqref{TheoremFirstInequalitySept2022}. Indeed they imply, together with Lemma \ref{lem:weakvsstrongNparticleMai2023}, that
\begin{align}
\notag \int_{\mathcal{P}_p(\mathcal{C}^d \times \mathcal{V})} \Gamma(\nu) dP(\nu) &\leq \liminf_{N \rightarrow +\infty} \E^{P_N} \left[ \Gamma(\hat{\nu}^N) \right] \\ 
\notag &\leq \liminf_{N \rightarrow +\infty} \mathcal{V}^N(t_0,\mathbf{x}_0^N) \\
&\leq \limsup_{N \rightarrow +\infty} \mathcal{V}^N(t_0,\mathbf{x}_0^N) \leq \mathcal{U}(t_0,\mu_0).
\label{11Mai202316h14}
\end{align}
However, by Lemma \ref{lem:weakvsstronglimitproblem} we have that $\mathcal{U}(t_0,\mu_0) = \inf_{m\in \mathcal{R}} \Gamma(m)$. Therefore, recalling that $P$ is supported on $\mathcal{R}$, we can deduce from \eqref{11Mai202316h14} that
$$ \int_{\mathcal{R}} \Gamma(\nu) dP(\nu) \leq \inf_{\nu \in \mathcal{R} } \Gamma(\nu),$$
which, in turn, implies that $P$ is supported on the set of optimal solutions for the limit problem and $\displaystyle \int_{\mathcal{P}_p(\mathcal{C}^d \times \mathcal{V})} \Gamma(\nu) dP(\nu) = \mathcal{U}(t_0,\mu_0)$. Finally, getting back to \eqref{11Mai202316h14}, we deduce that
$$\lim_{N \rightarrow +\infty} \mathcal{V}^N \bigl( t_0, \mu_0 \bigr) = \mathcal{U}(t_0,\mu_0),$$
which concludes the proof of the proposition.
\end{proof}

\section{Application to Large Deviations}

\label{sec: Application to Large Deviations Sept2022}

We are interested in the probability distribution of the first exit time from a given region of $\mathcal{P}_2(\R^d)$ for the empirical measure of a system of interacting particles. We assume that this region is described by $\Psi : \mathcal{P}_2(\R^d) \rightarrow \R$ as follows
$$ \Omega_{\infty} := \left \{ \mu \in \mathcal{P}_2(\R^d), \Psi( \mu ) < 0 \right \}, $$
and, for $(t,\mathbf{x}^N = (x^{1,N}, \dots, x^{N,N})) \in [0,T] \times  (\R^d)^N$ we introduce the probability
\begin{equation}
v^N(t,\mathbf{x}^N) := \mathbb{P} \left [ \forall s \in [0,t], \hat{\mu}_s^N \in \Omega_{\infty} \right ],
\label{defvN11Mai}
\end{equation}
where $(X_t^{1,N}, \dots, X_t^{N,N})$ is solution to the system of SDEs
\begin{equation*}
\left \{
\begin{array}{ll}
\displaystyle X_t^{i,N} = x^{i,N} + \int_0^t b(X_s^{i,N}, \hat{\mu}_s^{N}) ds +\sqrt{2}B_t^{i,N} &  t \in [0,T], \quad i \in \left \{1,\dots,N \right \}, \\
\displaystyle \hat{\mu}_s^{N} = \frac{1}{N} \sum_{i=1}^N \delta_{X_s^{i,N}}, 
\end{array}
\right.
\end{equation*}
with $(B_t^{1,N}, \dots, B_t^{N,N})$, $N$ independent $d$-dimensional standard Brownian motions supported on some probability space $(\Omega, \mathcal{F},  \mathbb{P})$.
 
The goal is to understand the asymptotic behavior of $v^N$ when $N \rightarrow +\infty$.

Throughout this section we take $L(x,q) = \frac{1}{2}|q|^2$ for all $(x,q) \in \R^d \times \R^d$ as well as $\mathcal{F} = \mathcal{G}=0$. We also make the following additional assumptions on the constraint.
\begin{assumption}
\begin{equation}
\mbox{The constraint  }\left \{ \Psi \leq 0 \right \}  \mbox{ is bounded in  }\mathcal{P}_1(\R^d).
\label{AdditionalBoundedConstraint}
\tag{APsibd}
\end{equation}
\end{assumption}
As a consequence, the constraints $\Omega_N := \left \{ (x_1, \dots, x_N) \in \R^{dN}, \Psi( \frac{1}{N} \sum_{i=1}^N \delta_{x_i} ) < 0 \right \}$ are bounded for all $N \geq 1$. We also assume that 
\begin{assumption}
\begin{equation}
\left \{
\begin{array}{ccc}
\displaystyle \mbox{For all $x \in \R^d$, $\displaystyle m \mapsto \frac{\delta \Psi}{\delta m}(m,x)$ is $\mathcal{C}^1$ with $\displaystyle (x,y) \mapsto \frac{\delta^2 \Psi}{\delta m^2}(m,x,y)$ } \\
\displaystyle \mbox{ in  $\mathcal{C}^2(\R^d \times \R^d)$ for all $m \in \mathcal{P}_2(\R^d)$ and $\displaystyle \frac{\delta^2 \Psi}{\delta m^2}(m,x,y)$ and its derivatives being} \\
\displaystyle \mbox{ jointly continuous and bounded in $\mathcal{P}_2(\R^d) \times \R^d \times \R^d$.}
\end{array}
\right.
\tag{APsiC2}
\label{APsiC2}
\end{equation}
\end{assumption}
and the transversality condition
\begin{assumption}
\begin{equation}
\begin{array}{ccc} 
\displaystyle
\int_{\R^d} |D_m\Psi(m,x) |^2 dm(x) \neq 0 \mbox{, whenever } \Psi(m)=0.
\end{array}
\tag{APsiTrans}
 \label{TransCondPsiMF2022}
\end{equation}
\end{assumption}
The condition \eqref{TransCondPsiMF2022} on the Wasserstein gradient of $\Psi$ at the boundary ensures that the closure $\overline{\Omega_N}$ of $\Omega_N$ in $(\R^d)^N$ is
$$\overline{\Omega_N} = \left \{ (x^1,\dots,x^N) \in (\R^d)^N, \Psi( \frac{1}{N} \sum_{i=1}^N \delta_{x^i} ) \leq 0 \right \}. $$
Similarly we have
$$ \overline{\Omega_{\infty}} = \left \{ \mu \in \mathcal{P}_2(\R^d), \Psi( \mu ) \leq 0 \right \}.$$

We will need further regularity for the constraint.
\begin{assumption}
\begin{align}
\notag &\displaystyle \frac{\delta^2 \Psi}{\delta m^2}  \mbox{ has a linear derivative, with bounded and jointly continuous first order} \\
&\mbox{ derivatives in the euclidean variables. }
\label{AdditionalC3}
\tag{APsiC3}
\end{align}
\end{assumption}
Under Assumption \eqref{AssumptionMeanFieldLimitSept2022} as well as these additional assumptions, for all $N \geq 1$, the constraint $\Omega_N$ is open, bounded and $\partial \Omega_N$ is a  manifold of class $\mathcal{C}^3$. These assumptions are satisfied for instance when $\displaystyle \Psi(m) = \int_{\R^d} \left(\sqrt{ |x-x_0|^2 +\delta^2 } - \delta \right) dm(x) - \kappa $ with $x_0 \in \R^d$, $\delta >0$ and $\kappa >0$.

Thanks to the additional assumptions \eqref{AdditionalBoundedConstraint}, \eqref{APsiC2}, \eqref{TransCondPsiMF2022} and \eqref{AdditionalC3} we are precisely in the framework of \cite{Fleming2006} section VI.6 -see also \cite{Patie2008}- and we can conclude that $v^N$, defined in \eqref{defvN11Mai} is $\mathcal{C}^{1,2}$ in $(0,T] \times \Omega_N$ and satisfies
\begin{equation*}
\left \{
\begin{array}{ll}
\displaystyle \partial_t v^N -\sum_{i=1}^N b^{i,N}(\mathbf{x}^N).D_{x^{i,N}}v^N - \Delta v^N = 0, & \mbox{ in } (0,T) \times \Omega_N \\
v^N(t,\mathbf{x}^N) = 0, & \mbox{ in } [0,T] \times \partial \Omega_N \\
v^N(0,\mathbf{x}^N) =1 & \mbox{ in } \Omega_N,
\end{array}
\right.
\end{equation*}
where $b^{i,N}(\mathbf{x}^N) = b(x^{i,N}, \frac{1}{N}\sum_{i=1}^N \delta_{x^{i,N}})$. Moreover, $v^N(t,\mathbf{x}) >0$ for all $(t,\mathbf{x}^N) \in (0,T] \times \Omega_N$. The connection with the control problems of the previous sections is given by the following proposition.

\begin{proposition} 
For all $\mathbf{x}_0^N \in \Omega_N$ and $t_0 \in [0,T]$ it holds
$$ \mathcal{V}^N(t_0, \mathbf{x}_0^N) = -\frac{2}{N} \log v^N(T-t_0,\mathbf{x}_0^N), $$
where $\mathcal{V}^N$ is the value function defined in Subsection \eqref{sec: The N-state problem} with $\mathcal{F}=\mathcal{G}=0$ and $L(x,q) = 1/2 |q|^2$.
\label{ControlRepresentation}
\end{proposition}

\begin{proof}
For all $(t_0, \mathbf{x}_0^N) \in [0,T) \times \Omega_N$ we define $w^N(t_0, \mathbf{x}_0^N) := -\frac{2}{N} \log v^N(T-t_0,\mathbf{x}_0^N)$. We are going to proceed by verification to show that $w^N(t_0, \mathbf{x}_0^N) = \mathcal{V}^N(t_0,\mathbf{x}_0^N)$ in $[0,T) \times \Omega_N$. For $\mathbf{x}_0^N=(x_0^{1,N}, \dots, x_0^{N,N}) \in \Omega_N$ and $t_0 \in [0,T)$, we define the following particle system where $\mathbf{X}_t^N = (X_t^{1,N}, \dots, X_t^{N,N})$
\begin{align*}
X_t^{i,N} &:= x_0^{i,N} - \int_{t_0}^{t \wedge \tau}N D_{x^{i}}w^N(t,\mathbf{X}_t^N) dt +\int_{t_0}^t b(X_s^{i,N}, \hat{\mu}^N_{\mathbf{X}_s^N})ds  + \sqrt{2} \int_{t_0}^{t \wedge \tau} d B^{i,N}_t \\
&= x_0^{i,N} + \int_{t_0}^{t \wedge \tau} 2\frac{D_{x^{i}}v^N(T-t,\mathbf{X}_t^N)}{v^N(T-t, \mathbf{X}_t^N)} dt +\int_{t_0}^t b(X_s^{i,N}, \hat{\mu}^N_{\mathbf{X}_s^N})ds+ \sqrt{2} \int_{t_0}^{t \wedge \tau} dB^{i,N}_t,
\end{align*}
where $\textbf{B}^N_t := ^t(B_t^{1,N}, \dots, B_t^{N,N})$ and $\tau$ is the first exit time from $\Omega_N$:
$$ \tau := \inf \{t \geq t_0, \mathbf{X}^N_t \notin \Omega_N \}. $$
For $\eta \geq 0$ small, we introduce the stopping time 
$$ \tau^{\eta} := \inf \{ t \geq t_0, v^N(T-t,\mathbf{X}^N_t) \leq \eta \}. $$
Notice that, by definition of $v^N$, it holds that $\tau^0 = \tau$.
Applying Itô's formula to $\log v^N(T-t, \mathbf{X}^N_t)$ yields, for $\eta >0$,
\begin{align*}
\log v^N(T- \tau^{\eta} &\wedge T, \mathbf{X}^N_{\tau^{\eta} \wedge T} ) = \log v^N(T,\mathbf{x}_0^N) \\
&+ \int_{t_0}^{\tau^{\eta} \wedge T} \left[ \frac{-\partial_t v^N}{v^N} + 2 \left|\frac{Dv^N}{v^N} \right|^2 + \sum_{i=1}^N \frac{D_{x^{i}}v^N}{v^N}.b^{i,N} + \frac{\Delta v^N}{v^N} - \left| \frac{Dv^N}{v^N} \right|^2 \right](T-t, \mathbf{X}_t^N) dt \\
&+ \int_{t_0}^{\tau^{\eta} \wedge T} \sqrt{2} \frac{Dv^N(T-t, \mathbf{X}^N_t)}{v^N(T-t, \mathbf{X}^N_t)}.d \textbf{B}_t^N \\
&=\log v^N(T,\mathbf{x}_0^N) + \int_{t_0}^{\tau^{\eta} \wedge T} \left| \frac{Dv^N(T-t, \mathbf{X}_t^N)}{v^N(T-t, \mathbf{X}^N_t)} \right|^2 dt \\
&+ \int_{t_0}^{\tau^{\eta} \wedge T} \sqrt{2} \frac{Dv^N(T-t, \mathbf{X}^N_t)}{v^N(T-t, \mathbf{X}^N_t)}.d \textbf{B}_t^N.
\end{align*}
Taking expectations and recalling the definition of $\tau^{\eta}$ we get
$$ \log (\eta) \mathbb{P}(\tau^{\eta} \leq T) + \mathbb{P}( \tau^{\eta} >T) \geq \log v^N(T, \mathbf{x}^N). $$
As a consequence, 
$$ \lim_{\eta \rightarrow 0} \mathbb{P}(\tau^{\eta} \leq T) =0 $$
and the control $-NDw^N(T-t, \mathbf{X}_t^N)$ is admissible. Let us show that it is optimal. Recalling the equation satisfied by $v^N$, it holds that
\begin{equation*}
\left \{
\begin{array}{ll}
-\partial_t w^N -\sum_{i=1}^N b^{i,N}(\mathbf{x}^N)D_{x^{i,N}}w^N+ \frac{N}{2} |Dw^N|^2 -\Delta w^N = 0, & \mbox{ in } (0,T) \times \Omega_N \\
w^N(t,\mathbf{x}^N)= +\infty, & \mbox{ in } [0,T] \times \partial \Omega_N \\
w^N(T,\mathbf{x}^N) =0 & \mbox{ in } \Omega_N.
\end{array}
\right.
\end{equation*}
Let us take another admissible control $\boldsymbol{\alpha}^N = (\alpha^{1,N},\dots, \alpha^{N,N})$ with the associated solution $\textbf{Y}^N = (Y_t^{1,N},\dots,Y_t^{N,N})$ to the SDE:
\begin{align*}
Y_t^{i,N} &:= x_0^{1,N} + \int_{t_0}^t \alpha^{i,N}_s ds + \int_{t_0}^t b^{i,N}(\textbf{Y}_s^N)ds  + \sqrt{2} \int_{t_0}^t dB^{i,N}_s.
\end{align*}
Being $\mathbf{\alpha}$ admissible, it holds that $\textbf{Y}_t^N$ belongs to $\Omega_N$ for all $t\in [t_0,T]$ almost surely. We can apply Itô's lemma to $w^N$ and get
\begin{align*}
0 = \E \left[w^N(T, \textbf{Y}_T^N) \right] &= w^N(t_0, \mathbf{x}_0^N) \\
&+  \E \left[ \int_{t_0}^T \partial_t w^N(t, \textbf{Y}_t^N) + \sum_{i=1}^N \left( \mathbf{\alpha}^{i,N}_t +b^{i,N}(\textbf{Y}_t^N) \right) .D_{x^{i}}w^N(t, \textbf{Y}_t^N) +  \Delta w^N(t, \textbf{Y}_t^N) \right] dt \\
& = w^N(t_0, \mathbf{x}_0^N) + \E \left[ \int_{t_0}^T\left( \boldsymbol{\alpha}^N_t . Dw^N(t, \textbf{Y}_t^N) + \frac{N}{2} |Dw^N(t, \textbf{Y}_t^N)|^2 \right) dt \right] \\
& \geq w^N(t_0, \mathbf{x}^N) - \E \left[ \int_{t_0}^T \frac{1}{2N} |\boldsymbol{\alpha}^N_t |^2 dt \right]
\end{align*}
with equality if and only if $\boldsymbol{\alpha}^N_t = - NDw^N(t,\textbf{Y}^N_t)$. 
This means that the control $- NDw^N(t,\textbf{Y}^N_t)$ is optimal and the optimal value is given by $w^N(t_0, \mathbf{x}_0^N)$ which concludes the proof of the proposition.
\end{proof}

Notice that a by-product of Proposition \ref{ControlRepresentation} is to characterize $\mathcal{V}^N$ as the unique solution in $\mathcal{C}^{1,2}([0,T) \times \Omega_N)$ to the HJB equation
\begin{equation*}
\left \{
\begin{array}{ll}
-\partial_t \mathcal{V}^N -\sum_{i=1}^N b^{i,N}(\mathbf{x}^N).D_{x^{i,N}}\mathcal{V}^N+ \frac{N}{2}\sum_{i=1}^N |D_{x^{i,N}}\mathcal{V}^N|^2 -\sum_{i=1}^N\Delta_{x^{i,N}} \mathcal{V}^N = 0, & \mbox{ in } (0,T) \times \Omega_N \\
\mathcal{V}^N(t,\mathbf{x}^N) = +\infty, & \mbox{ in } [0,T] \times \partial \Omega_N \\
\mathcal{V}^N(T,\mathbf{x}^N) =0 & \mbox{in } \Omega_N.
\end{array}
\right.
\end{equation*}
The same argument extends without difficulty when additional mean-field costs $\mathcal{F}$ and $\mathcal{G}$ satisfying Assumption \eqref{Ureg} are considered.

Combining Proposition \ref{ControlRepresentation} and Theorem \ref{MainTheoremMFconstraint2022} we obtain the following convergence.

\begin{corollary}
Let Assumption \eqref{AssumptionMeanFieldLimitSept2022} as well as Assumptions \eqref{AdditionalBoundedConstraint}, \eqref{APsiC2}, \eqref{TransCondPsiMF2022} and \eqref{AdditionalC3} hold. Assume that $\Psi(\mu_0) <0$ and write $\mathbf{x}_0^N = (x^{1,N}_0,\dots, x^{N,N}_0)$. Then it holds
$$ \lim_{N \rightarrow + \infty} \frac{2}{N} \log v^N(T-t_0, \mathbf{x}_0^N) = - \mathcal{U}(t_0,\mu_0),$$
whenever $\displaystyle  \frac{1}{N} \sum_{i=1}^N \delta_{x^{i,N}_0} \rightarrow \mu_0$ is $\mathcal{P}_2(\R^d)$. 
\label{Coro1_10MAI}
\end{corollary}

This is a special case of the general result of Dawson of Gärtner, \cite{Dawson1987}. Contrary to the general result of \cite{Dawson1987}, we don't need to express the Large Deviation principle with $\limsup$ and $\liminf$. This is due to the fact that the constraint is “regular" with respect to the rate function, as can be seen with the stability result of Section \ref{sec: Stability section}. The optimality conditions of Theorem \ref{TheoremOptimalityConditions} give a new way to compute the limit $-\mathcal{U}(t_0,\mu_0)$. 

\begin{corollary}
Under the assumptions of Corollary \eqref{Coro1_10MAI} we have
$$ \lim_{N \rightarrow + \infty} \frac{2}{N} \log v^N \bigl(T-t_0, \mathbf{x}_0^N \bigr) = - \int_{\R^d} u(t_0,x)d\mu_0(x),$$
for a solution $(u,\mu,\nu,\eta)$ of the optimality conditions
\begin{equation}
\left\{
\begin{array}{lr}
\displaystyle -\partial_t u(t,x) + \frac{1}{2}|Du(t,x)|^2 - Du(t,x).b \bigl(x,\mu(t) \bigr) - \Delta u(t,x) \\
\displaystyle \hspace{10pt} = \nu(t) \frac{\delta \Psi}{\delta m}\bigl(\mu(t),x \bigr) + \int_{\R^d} Du(t,y). \frac{\delta b}{\delta m}\bigl(y, \mu(t),x \bigr)d\mu(t)(y) & \mbox{in   } (t_0,T)\times \R^d, \\
\displaystyle \partial_t \mu - \mdiv \bigl( Du(t,x) \mu \bigr) +\mdiv \bigl(b(x,\mu(t) )\mu \bigr)- \Delta \mu = 0  & \mbox{in   } (t_0,T)\times \R^d,  \\
\displaystyle \mu(t_0) = \mu_0, \hspace{30pt} \displaystyle u(T,x) = \eta \int_{\R^d} \frac{\delta \Psi}{\delta m} \bigl(\mu(T),x \bigr)  & \mbox{       in      } \R^d.
\end{array}
\right.
\end{equation}
In the above optimality conditions $\mu$ belongs to $\mathcal{C}^0([t_0,T], \mathcal{P}_2(\R^d))$, $\nu \geq 0$ belongs to $L^{\infty}([t_0,T])$, $\beta$ to $\R^+$ and finally $u$ belongs to $\mathcal{C}([t_0,T], E_n)$ and $Du$ is bounded and globally Lipschitz continuous. Moreover $\lambda$ and $\beta$ satisfy the exclusion conditions
$$ \int_{t_0}^T \Psi(\mu(t)) \nu(t)dt =0, \quad \quad \eta \Psi(\mu(T)) =0.$$
\end{corollary}

\begin{proof}
The optimality conditions are given by Proposition \ref{OptimalityConditionMainTheorem2021} with a Lagrange multiplier $\nu' \in \mathcal{M}^+([t_0,T])$. Thanks to Assumptions \eqref{APsiC2} and \eqref{TransCondPsiMF2022} we can apply Theorem 2.2 of \cite{Daudin2021} (which applies equally well when $b$ satisfies Assumption \eqref{Ab})  to infer that the Lagrange multiplier $\nu'$ has the form $\nu' = \nu + \eta \delta_T$ for some $\nu \in L^{\infty}([t_0,T])$ and some $\eta \in \R^+$. The fact that $\eta$ appears in the terminal condition for $u$ follows then from the representation formula \eqref{RepresentationformulaHJB10Mai}. Finally the additional regularity of $u$ and $Du$ follows from the boundness of $\nu$ and \eqref{RepresentationformulaHJB10Mai} again, see \cite{Daudin2021} Theorem 1.1.
\end{proof}

\section{Appendix}

\subsection{Optimality conditions}

\label{subsec:OptimalityConditions}

We will need some preliminary facts.

\begin{lemma}
There is an admissible couple $(\bar{\mu},\bar{\alpha})$ satisfying \eqref{FPwithdriftb10janv2023} and $J(t_0,\mu_0,(\bar{\mu},\bar{\alpha}))<+\infty$ such that $\Psi(\bar{\mu}(t)) \leq -\eta $ for some $\eta >0$ and all $t \in [t_0,T]$. 
\label{controllabilitylemma30Mars2023}
\end{lemma}

\begin{proof}
Following \cite{Daudin2021} Lemma 4.1, for every $\epsilon >0$, we can build a solution $(\bar{\mu},\bar{\beta})$ to 
$$ \partial_t \mu + \mdiv(\beta \mu) - \Delta \mu = 0 \hspace{15pt}  \mbox{ in } (t_0,T) \times \R^d, \hspace{30pt} \mu(t_0)=\mu_0, $$
such that $d_2(\bar{\mu}(t), \mu_0) \leq \epsilon$ for all $t \in [t_0,T]$ and $\int_{t_0}^T\int_{\R^d} |\bar{\beta}(t,x)|^2 d\bar{\mu}(t)(x)dt <+\infty$. In particular, since $\Psi(\mu_0) <0$, we can take $\epsilon$ small enough so that $\Psi(\bar{\mu}(t)) \leq -\Psi(\mu_0)/2$ for all $t \in [t_0,T]$. Thanks to the growth assumption on $L$, and the boundness of $b$, we find that $\bigl(\bar{\mu},\bar{\beta}(t,x) - b(x,\bar{\mu}(t)) \bigr)$ satisfies the desired properties. 
\end{proof}

\begin{lemma}
There exists a solution $(\tilde{\mu},\tilde{\alpha})$ to Problem \eqref{ConstrainedProblem}.
\end{lemma}

\begin{proof}
Take $(\mu_n,\alpha_n)_{n \geq 0}$ a minimizing sequence such that $J \bigl(t_0,\mu_0;(\mu_n,\alpha_n) \bigr) \leq \inf_{(\mu,\alpha)} J \bigl(t_0,\mu_0,(\mu,\alpha) \bigr) +1$ for all $n \in \mathbb{N}$. Using the previous lemma, the growth condition on $L$, the boundness of $b$ and that $\mathcal{F}$ and $\mathcal{G}$ are bounded from below, we find that 
\begin{equation} 
\int_{t_0}^T \int_{\R^d}  |\alpha_n(t,x) |^2d \mu_n(t)(x)dt \leq C
\label{energybound30Mars2023}
\end{equation}
for some $C>0$ and all $n \geq 0$. By classical argument, see Proposition 2.1 in \cite{Daudin2021}, this implies that 
$$ \sup_{t \in [t_0,T]} \int_{\R^d} |x|^2d\mu_n(t) + \sup_{t  \neq s \in [t_0,T]} \frac{d_2(\mu_n(t),\mu_n(s))}{\sqrt{|t-s|}} \leq C$$
for some $C >0$ independent from $n$. From \eqref{energybound30Mars2023} we also deduce, by Cauchy-Schwartz inequality, that we have
\begin{align*} 
\int_{t_0}^T \int_{\R^d}  |\alpha_n(t,x)|\mu_n(t)(x)dt &\leq \sqrt{T-t_0} \sqrt{ \int_{t_0}^T \int_{\R^d} |\alpha_n(t,x) |^2d\mu_n(t)(x)dt } \\
&\leq \sqrt{TC}
\end{align*}
and therefore $\alpha_n\mu_n$ is bounded in $\mathcal{M}([t_0,T] \times \R^d,\R^d)$ uniformly in $n \in \mathbb{N}$. Now we take $\delta \in (0,1)$ and we apply Banach-Alaoglu and Ascoli theorems to deduce that, up to taking a sub-sequence, $(\mu_n, \mu_n \alpha_n)_{n \geq 0}$ converges in $\mathcal{C}([t_0,T] , \mathcal{P}_{2-\delta}(\R^d)) \times \mathcal{M}([t_0,T] \times \R^d,\R^d)$ toward some $(\tilde{\mu},\tilde{\omega}) \in \mathcal{C} \bigl([t_0,T],\mathcal{P}_2(\R^d) \bigr) \times \mathcal{M} \bigl([t_0,T] \times \R^d,\R^d \bigr)$. Using Theorem 2.34 and Example 2.36 in \cite{Ambrosio2000} to handle the term involving the Lagrangian $L$, we conclude that $\tilde{\omega}$ is absolutely continuous with respect to $\tilde{\mu}(t) dt$ and, taking $\tilde{\alpha} := \frac{d \tilde{\omega}}{dt \otimes d\tilde{\mu}}$, that 
$$ J \bigl(t_0,\mu_0;(\tilde{\mu},\tilde{\alpha}) \bigr) \leq \liminf_{n \rightarrow + \infty} J \bigl(t_0,\mu_0,(\mu_n,\alpha_n) \bigr). $$
We easily check that $(\tilde{\mu},\tilde{\alpha})$ satisfies the Fokker-Planck equation and the state constraint and is therefore a solution to Problem \eqref{ConstrainedProblem}.
\end{proof}

\begin{lemma}
If $(\tilde{\mu},\tilde{\alpha})$ is a solution to Problem \eqref{ConstrainedProblem}, then $(\tilde{\mu},\tilde{\beta}):= (\tilde{\mu},\tilde{\alpha} + b(x, \tilde{\mu}(t)))$ is a solution to 
\begin{equation}
 \inf_{(\mu,\beta)} J'^{l} \bigl(t_0,\mu_0; (\mu,\beta) \bigr)
\end{equation}
with $J'^{l} \bigl(t_0,\mu_0;(\mu,\beta) \bigr)$ defined by
\begin{align*}
 J'^l \bigl(t_0,\mu_0;(\mu,\beta) \bigr) &:= \int_{t_0}^T \int_{\R^d} L\Bigl(x,\beta(t,x)-b(\tilde{\mu}(t),x) \Bigr)d\mu(t)(x) dt \\
 &- \int_{t_0}^T \int_{\R^d} \int_{\R^d} \partial_q L \bigl(y,\tilde{\alpha}(t,y) \bigr).\frac{\delta b}{\delta m} \bigl(\tilde{\mu}(t),y,x \bigr)d\tilde{\mu}(t)(y)d\mu(t)(x)dt \\
&+ \int_{t_0}^T \int_{\R^d} \frac{\delta \mathcal{F}}{\delta m}\bigl (\tilde{\mu}(t),x \bigr)d\mu(t)(x)dt  + \int_{\R^d} \frac{\delta \mathcal{G}}{\delta m} \bigl(\tilde{\mu}(T),x \bigr)d\mu(T)(x)
\end{align*}
where the infimum is taken over the couples $(\mu, \beta) \in \mathcal{C}([t_0, T], \mathcal{P}_2(\R^d)) \times L^2_{dt \otimes \mu(t)} ([t_0,T] \times \R^d, \R^d)$ satisfying in the sense of distributions the Fokker-Planck equation 
\begin{equation}
\left \{
\begin{array}{ll}
 \partial_t \mu + \mdiv \bigl(\beta(t,x) \mu \bigr) - \Delta \mu = 0 & \mbox{ in }(t_0,T) \times \R^d \\
\mu(t_0)=\mu_0,
\end{array}
\right.
\label{FPEforbeta30Mars2023}
\end{equation}
under the constraint that $\Psi(\mu(t)) \leq 0$ for all $t \in [t_0,T]$.
\label{lemmaLinearization30mars2023}
\end{lemma}

\begin{proof}
We first recall that, being $(\tilde{\mu},\tilde{\alpha})$ a solution to Problem \eqref{ConstrainedProblem}, then $(\tilde{\mu},\tilde{\beta})$ is a solution to Problem \eqref{ConstrainedProblemFormulationbeta}.  Now we take $(\mu,\beta)$ satisfying  \eqref{FPEforbeta30Mars2023} and such that  $J'\bigl(t_0,\mu_0,(\mu,\beta) \bigr) < +\infty$. We define $\omega = \beta \mu$ and $\tilde{\omega} = \tilde{\beta} \tilde{\mu}$. We take $\lambda \in (0,1)$ and we let  $ (\omega_\lambda, \mu_\lambda) = (1-\lambda) (\tilde{\omega}, \tilde{\mu}) + \lambda (\omega, \mu)$. In particular, $\omega_{\lambda}$ is absolutely continuous with respect to $\mu_{\lambda}(t)\otimes dt$ and we let $\beta_{\lambda} = \frac{d\omega_{\lambda}}{dt \otimes \mu_{\lambda}(t)}$.  By convexity of the constraint and linearity of the Fokker-Planck equation \eqref{FPEforbeta30Mars2023} in $(\mu, \beta \mu)$, the couple $(\mu_{\lambda},\beta_{\lambda})$ satisfies \eqref{FPEforbeta30Mars2023} as well as the state constraint. By minimality of $(\tilde{\mu},\tilde{\beta})$ it holds that
$$ J' \bigl(t_0,\mu_0,(\mu_{\lambda},\beta_{\lambda}) \bigr) \geq J' \bigl(t_0,\mu_0, (\tilde{\mu},\tilde{\beta}) \bigr). $$
On the other hand, by convexity of $ \R^d \times \R^+ \ni (\beta ,m) \mapsto L\Bigl(x,\frac{\beta}{m} \Bigr)m $ (set to be $+\infty$ if $m =0$) for all $x \in \R^d$, we have
\begin{align*} 
J' \bigl(t_0,\mu_0,&(\mu_{\lambda},\beta_{\lambda}) \bigr) \leq \int_{t_0}^T \int_{\R^d} L \Bigl(x, \tilde{\beta}(t,x)- b \bigl(\mu_{\lambda}(t),x \bigr) \Bigr)d\tilde{\mu}(t)(x)dt \\
&+ \lambda \int_{t_0}^T \int_{\R^d} L\Bigl(x, \beta(t,x)-b \bigl(\mu_{\lambda}(t),x \bigr) \Bigr)d\mu(t)(x)dt \\
&-\lambda \int_{t_0}^T \int_{\R^d} L \Bigl(x, \tilde{\beta}(t,x)-b \bigl(\mu_{\lambda}(t),x \bigr) \Bigr)d\mu(t)(x)dt  \\
&+ \int_{t_0}^T  \mathcal{F}\bigl(\mu_{\lambda}(t) \bigr)dt + \mathcal{G} \bigl(\mu_{\lambda}(T) \bigl). 
\end{align*}
Combining the two inequalities we find that, for all $\lambda \in (0,1)$,
\begin{align*}
 &\int_{t_0}^T \int_{\R^d} L \Bigl(x, \beta(t,x)-b \bigl(\mu_{\lambda}(t),x \bigr) \Bigr)d\mu(t)(x)dt -\int_{t_0}^T \int_{\R^d} L \Bigl(x, \tilde{\beta}(t,x)-b \bigl(\mu_{\lambda}(t),x \bigr) \Bigr)d\mu(t)(x)dt \\
 & \geq \frac{1}{\lambda} \left[ \int_{t_0}^T \int_{\R^d} L \Bigl(x, \tilde{\beta}(t,x)- b \bigl(\tilde{\mu}(t),x \bigr) \Bigr)d\tilde{\mu}(t)(x)dt - \int_{t_0}^T \int_{\R^d} L \Bigl(x, \tilde{\beta}(t,x)- b \bigl(\mu_{\lambda}(t),x \bigr) \Bigr)d\tilde{\mu}(t)(x)dt  \right] \\
 & + \int_{t_0}^T \frac{1}{\lambda} \left[ \mathcal{F} \bigl(\tilde{\mu}(t) \bigr) - \mathcal{F} \bigl(\mu_{\lambda}(t) \bigr) \right]dt +  \frac{1}{\lambda} \left[ \mathcal{G} \bigl(\tilde{\mu}(T) \bigr) - \mathcal{G} \bigl(\mu_{\lambda}(T) \bigr) \right].
 \end{align*}
Letting $\lambda$ tends to $0$, using the differentiability of the mean-field costs and rearranging the terms gives
$$ J'^l \bigl( t_0,\mu_0; (\tilde{\mu}, \tilde{\beta}) \bigr) \leq  J'^l \bigl( t_0,\mu_0; (\mu, \beta) \bigr)$$
which concludes the proof of the Lemma.
\end{proof}

The next step to derive the optimality conditions relies on the following form of Von Neumann min/max theorem, see \cite{Simons1998} for the proof under additional compactness assumptions and the appendix of \cite{Orrieri2019} for this version.
\begin{theorem}(Von Neumann)
\label{VN}
Let $\mathbb{A}$ and $\mathbb{B}$ be convex sets of some vector spaces and suppose that $\mathbb{B}$ is endowed with some Hausdorff topology. Let $\cL$ be a function satisfying :
$$ a \rightarrow \cL (a,b) \mbox{      is concave in }\mathbb{A} \mbox{ for every }b \in \mathbb{B}, $$
$$ b \rightarrow \cL (a,b) \mbox{      is convex in }\mathbb{B} \mbox{ for every }a \in \mathbb{A}. $$
Suppose also that there exists $a_* \in \mathbb{A}$ and $C_* > \sup_{a \in \mathbb{A} } \inf_{b \in \mathbb{B}} \cL(a,b)$ such that :
$$ \mathbb{B}_*:= \left \{ b \in \mathbb{B}, \cL (a_*,b) \leq C_* \right \} \mbox{      is not empty and compact in } \mathbb{B}, $$
$$ b \rightarrow \cL(a,b) \mbox{        is lower semicontinuous in } \mathbb{B}_* \mbox{   for every  } a \in \mathbb{A}. $$
Then,
$$ \min_{b\in \mathbb{B} } \sup_{a \in \mathbb{A} } \cL(a,b) = \sup_{ a \in \mathbb{A} } \inf_{b \in \mathbb{B} } \cL(a,b). $$
\end{theorem}

\begin{lemma}
If $(\tilde{\mu},\tilde{\alpha})$ is a solution to Problem \eqref{ConstrainedProblem}, then there exists some $\tilde{\nu} 
\in \mathcal{M}^+([t_0,T])$ such that 
\begin{equation}
\Psi(\tilde{\mu}(t)) = 0 \hspace{ 30pt} \tilde{\nu} -ae
\label{Exclusion12MAi2023}
\end{equation}
and $(\tilde{\mu},\tilde{\beta}):= \bigl(\tilde{\mu},\tilde{\alpha} + b(x, \tilde{\mu}(t)) \bigr)$ is a solution to 
\begin{equation} 
\inf_{(\mu,\beta)} J'^l \bigl(t_0,\mu_0,(\mu,\beta) \bigr) + \int_{t_0}^T \Psi \bigl(\mu(t) \bigr)d \tilde{\nu}(t)
\label{Toprove12Mai2023}
\end{equation}
where $J'^l$  was defined in Lemma \ref{lemmaLinearization30mars2023} and the infimum is taken over the couples $(\mu,\beta)$ satisfying \eqref{FPEforbeta30Mars2023} but not necessarily the state constraint.
\label{lem:minmax11Mai2023}
\end{lemma}

\begin{proof}
We set-up the min-max argument. We define $\mathbb{A}$ as the subset of $\mathcal{C}([t_0,T], \mathcal{P}_2(\R^d)) \times \mathcal{M}([t_0,T] \times \R^d,\R^d)$ consisting of elements $(\mu, \omega)$ such that
\begin{equation}
\left\{
\begin{array}{ll}
\displaystyle \omega \mbox{ is absolutely continuous w.r.t. } \mu(t) \otimes dt, \quad \frac{d\omega}{dt \otimes \mu(t)} \in L^2_{dt \otimes d\mu(t)} \bigl([t_0,T] \times \R^d, \R^d \bigr), \\
\partial_t \mu +\mdiv(\omega) - \Delta \mu = 0\hspace{15pt} \mbox{in } (t_0,T) \times \R^d ,\\
\mu(t_0)=\mu_0,
\end{array}
\right.
\end{equation}
where the Fokker-Planck equation is understood in the sense of distributions. We also set $\mathbb{B} = \mathcal{M}^+([t_0,T])$. We define $\mathcal{L} : \mathbb{A} \times \mathbb{B} \rightarrow \R$ by
$$ \mathcal{L} \bigl ((\mu,\omega), \nu \bigr) = J'^l \bigl(t_0,\mu_0,(\mu,\omega) \bigr) + \int_{t_0}^T \Psi \bigl(\mu(t) \bigr)d\nu(t), $$
where we set, by abuse of notation, $J'^l\bigl(t_0,\mu_0,(\mu,\omega) \bigr) = J'^l\Bigl(t_0,\mu_0, \bigl(\frac{\omega}{dt \otimes d\mu},\mu\bigr) \Bigr)$. It is clear that 
\begin{equation}
\inf_{(\mu,\beta)} J'^l \bigl(t_0,\mu_0,(\mu,\beta) \bigr) = \inf_{(\mu,\omega) \in \mathbb{A}} \sup_{\nu \in \mathbb{B}} \mathcal{L} \bigl((\mu,\omega),\nu \bigr)
\label{Troptard30Mars2023}
\end{equation}
where the first infimum is taken over the couples $(\mu,\beta)$ satisfying \eqref{FPEforbeta30Mars2023} as well as the state constraint. It is plain to check that, for every $\nu \in \mathbb{B}$, $\mathcal{L}(., \nu)$ is convex over $\mathbb{B}$ and that, for any $(\mu,\omega) \in \mathbb{B}$, $\mathcal{L} \bigl((\mu,\omega),. \bigr)$ is concave (linear in fact) and continuous over $\mathbb{A}$. Moreover, arguing as in  Lemma  \ref{controllabilitylemma30Mars2023}, we find $\bigl(\bar{\mu},\bar{\omega} \bigr) \in \mathbb{A}$ such that $J'^l \bigl(t_0,\mu_0,(\bar{\mu}, \bar{\omega}) \bigr) < +\infty$ and $\Psi \bigl(\bar{\mu}(t) \bigr) \leq -\eta$ for some $\eta>0$ and all $t \in [t_0,T]$. As a consequence, using the continuity of $\nu \mapsto \mathcal{L}\bigl( (\bar{\mu}, \bar{\omega}), \nu \bigr)$, we find that the set 
\begin{align*}
\Bigl \{ \nu \in \mathcal{M}^+&([t_0,T]),  \quad \mathcal{L}\bigl((\bar{\mu},\bar{\omega}),\nu \bigl) \geq     \inf_{(\mu,\omega) \in \mathbb{A}} \sup_{\nu \in \mathbb{B}} \mathcal{L}\bigl((\mu,\omega),\nu \bigr) \Bigr \} \\
&\subset \left \{  \nu \in \mathcal{M}^+([t_0,T]), \quad \nu([t_0,T]) \leq \frac{J'^l \bigl( t_0, \mu_0; (\bar{\mu},\bar{\omega}) \bigr) -  \inf_{(\mu,\omega) \in \mathbb{A}} \sup_{\nu \in \mathbb{B}} \mathcal{L} \bigl((\mu,\omega),\nu \bigr) }{\eta } \right \} 
\end{align*}
is non-empty and compact in $\mathcal{M}^+([t_0,T])$. Moreover, for all $(\mu,\omega) \in \mathbb{B}$, $\nu \mapsto \mathcal{L}\bigl((\mu,\omega), \nu \bigr)$ is continuous over $\mathbb{B}$. Applying Von-Neumann min-max theorem, we find that 
\begin{equation} 
\inf_{(\mu,\omega) \in \mathbb{B}} \sup_{\nu \in \mathbb{A}} \mathcal{L} \bigl((\mu,\omega),\nu \bigr) = \max_{\nu \in \mathbb{A}} \inf_{(\mu,\omega) \in \mathbb{B}} \mathcal{L} \bigl((\mu,\omega),\nu \bigr).
\label{MinMax12MAi2023}
\end{equation}
Let $\tilde{\nu}$ be a solution to the dual problem, ie $\max_{\nu \in \mathbb{A}} \inf_{(\mu,\omega) \in \mathbb{B}} \mathcal{L}\bigl((\mu,\omega),\nu \bigr) = \inf_{(\mu,\omega) \in \mathbb{B}} \mathcal{L} \bigl((\mu,\omega),\tilde{\nu} \bigr)$. Combining \eqref{Troptard30Mars2023} and \eqref{MinMax12MAi2023} we deduce that $(\tilde{\mu}, \tilde{\beta})$ is a solution to \eqref{Toprove12Mai2023}. It remains to prove \eqref{Exclusion12MAi2023}. By \eqref{Troptard30Mars2023} it holds that
$$J'^l \bigl(t_0,\mu_0; (\tilde{\mu},\tilde{\beta}) \bigr) = \inf_{(\mu,\omega) \in \mathbb{B}} \mathcal{L} \bigl((\mu,\omega),\tilde{\nu} \bigr) \leq J'^l \bigl(t_0,\mu_0, (\tilde{\mu},\tilde{\beta}) \bigr) + \int_{t_0}^T \Psi \bigl(\tilde{\mu}(t) \bigr)d\tilde{\nu}(t). $$
This implies that $\int_{t_0}^T \Psi \bigl(\tilde{\mu}(t) \bigr)d\tilde{\nu}(t)  \geq 0$ but $\Psi \bigl(\tilde{\mu}(t) \bigr) \leq 0$ for all $t\in [t_0,T]$ and therefore,
$$ \Psi \bigl(\tilde{\mu}(t) \bigr) = 0 \hspace{30pt} \tilde{\nu}-ae \mbox{ in } [t_0,T], $$
which concludes the proof of the lemma.
\end{proof}
We can finally prove the optimality conditions of Proposition \ref{TheoremOptimalityConditions}.

\begin{proof}[Proof of Proposition \ref{TheoremOptimalityConditions}]

Take $(\tilde{\mu}, \tilde{\alpha})$ a solution to Problem \eqref{ConstrainedProblem}. Using  Lemma \ref{lem:minmax11Mai2023} and arguing as in Lemma \ref{lemmaLinearization30mars2023}, we find that $\bigl(\tilde{\mu},\tilde{\beta} \bigr) := \bigl(\tilde{\mu},\tilde{\alpha} - b(x,\tilde{\mu}(t)) \bigr)$ is solution to 
$$ \inf_{(\mu,\beta)} J'^l \bigl(t_0,\mu_0;(\mu,\beta) \bigr) + \int_{t_0}^T\int_{\R^d} \frac{\delta \Psi}{\delta m} \bigl(\tilde{\mu}(t),x \bigr)d\mu(t)(x) d\nu(t), $$
where the infimum is taken over the couples $(\mu,\beta)$ satisfying \eqref{FPEforbeta30Mars2023} but not necessarily the state constraint. This means that $(\tilde{\mu},\tilde{\alpha})$ is solution to 
$$ \inf_{(\mu,\alpha)} J^l(t_0,\mu_0,(\mu,\alpha)) $$
where $J^{l}$ is defined by
\begin{align*}
 J^l \bigl(t_0,\mu_0;(\mu,\alpha) \bigr) &:= \int_{t_0}^T \int_{\R^d} L \bigl(x,\alpha(t,x) \bigr)d\mu(t)(x) dt \\
 &- \int_{t_0}^T \int_{\R^d} \int_{\R^d} \partial_q L \bigl(y,\tilde{\alpha}(t,y) \bigr).\frac{\delta b}{\delta m} \bigl(\tilde{\mu}(t),y,x \bigr)d\tilde{\mu}(t)(y)d\mu(t)(x)dt \\
&+ \int_{t_0}^T \int_{\R^d} \frac{\delta \mathcal{F}}{\delta m}\bigl(\tilde{\mu}(t),x \bigr)d\mu(t)(x)dt  + \int_{\R^d} \frac{\delta \mathcal{G}}{\delta m} \bigl(\tilde{\mu}(T),x \bigr)d\mu(T)(x) \\
&+ \int_{t_0}^T \int_{\R^d} \frac{\delta \Psi}{\delta m} \bigl(\tilde{\mu}(t),x \bigr) d\mu(t)(x)d\nu(t)
\end{align*}
and the infimum is taken over the couples $(\mu, \alpha) \in \mathcal{C}([t_0, T], \mathcal{P}_2(\R^d)) \times L^2_{dt \otimes \mu(t)} ([t_0,T] \times \R^d, \R^d)$ satisfying in the sense of distributions the Fokker-Planck equation 
\begin{equation}
\left \{
\begin{array}{ll}
 \partial_t \mu + \mdiv(\alpha \mu) +  \mdiv \bigl(b(\tilde{\mu}(t),x) \mu \bigr) - \Delta \mu = 0 & \mbox{ in }(t_0,T) \times \R^d, \\
\mu(t_0)=\mu_0.
\end{array}
\right.
\end{equation}
We are now dealing with a standard control problem (except maybe for the presence of the measure $\nu$) for a linear Fokker-Planck equation and, importantly, without state constraint. Therefore, arguing by verification --- see e.g. \cite{Daudin2021} proof of Theorem 2.3 for the detailed computation---, we find that $\tilde{\alpha}= -\partial_pH(x,Du)$ where $u \in L^{\infty}([t_0,T], E_n)$ is solution in the sense of Definition \eqref{definitionHJBsol2023} to 
\begin{equation*}
\left \{
\begin{array}{ll}
\displaystyle -\partial_t u -\Delta u + H(x,Du) - b \bigl(\tilde{\mu}(t),x \bigr).Du =  \nu(t)  \frac{\delta \Psi}{\delta m} \bigl(\tilde{\mu}(t),x \bigr) \\
\displaystyle \hspace{10pt}  - \int_{\R^d} \partial_q L \bigl (y,\tilde{\alpha}(t,y) \bigr).\frac{\delta b}{\delta m} \bigl(\tilde{\mu}(t),y,x \bigr)d\tilde{\mu}(t)(y) +  \frac{\delta \mathcal{F}}{\delta m}\bigl(\tilde{\mu}(t),x \bigr) & \mbox{ in } (t_0,T) \times \R^d, \\
\displaystyle u(T,x) = \frac{\delta \mathcal{G}}{\delta m}\bigl (\tilde{\mu}(T),x \bigr) & \mbox{in } \R^d.
\end{array}
\right.
\end{equation*}
The existence of such a solution $u$ is guaranteed by Theorem 5.1. in \cite{Daudin2021}. Collecting the equations satisfied by $u$ and $\mu$ as well as noticing that 
$$\partial_q L \bigl(x,\tilde{\alpha}(t,x) \bigr) = \partial_qL \Bigl(x, -\partial_pH\bigl(x,Du(t,x) \bigr) \Bigr) = -Du(t,x)$$
gives the optimality conditions.
 \end{proof} 

\subsection{Concentration limit}

\label{sec: Concentration limit}

We consider $\mathbf{x}_0 = (x_0^1,\dots,x_0^N) \in (\R^d)^N$ such that $\frac{1}{N}\sum_{i=1}^N \delta_{x_0^{i}} \rightarrow \mu_0$ in $\mathcal{P}_2(\R^d)$ as $N \rightarrow +\infty$. For $b : [0,T] \times \R^d \times \mathcal{P}_1(\R^d) \rightarrow \R^d$ uniformly Lipschitz continuous we consider the particle system
\begin{equation*}
\left \{
\begin{array}{ll}
\displaystyle dX_t^{i} = b(X_t^{i}, \hat{\mu}_t^N)dt + \sqrt{2}dB_t^{i} & 1\leq i \leq N, t\in [0,T] \\
\displaystyle \hat{\mu}_t^N = \frac{1}{N}\sum_{i=1}^N \delta_{X_t^{i}} \\
\displaystyle X_0^{i} = x_0^{i}, & 1\leq i \leq N
\end{array}
\right.
\end{equation*}
where $B^1,\dots,B^N$ are independent Brownian motions.

We also consider $\mu \in \mathcal{C}([0,T], \mathcal{P}_2(\R^d))$  solution to 
\begin{equation}
\left \{
\begin{array}{ll}
\displaystyle \partial_t \mu +\mdiv(b(x,\mu(t))\mu) - \Delta \mu = 0 &\mbox{ in } (0,T) \times \R^d \\
\mu(0) = \mu_0.
\end{array}
\right.
\label{FPeAppendix2023}
\end{equation}

\begin{proposition} In this setting, it holds
\begin{equation} 
\lim_{N \rightarrow +\infty} \E \left[\sup_{t \in [0,T]} d_1( \hat{\mu}_t^N, \mu(t)) \right] =0.
\label{ToproveintheAppendixConcentration}
\end{equation}
\label{PropconcentrationAppendix}
\end{proposition}

\begin{proof}
By classical arguments, see Oelschlager \cite{Oelschlager1984}, $(\mathcal{L}(\hat{\mu}^N_.))_{n \in \mathbb{N}}$ is pre-compact in $\mathcal{P}( \mathcal{C}([0,T], \mathcal{P}_1(\R^d))$. The limit points are supported on the set of solutions to the Fokker-Planck equation \eqref{FPeAppendix2023}. This equation admits a unique solution in $\mathcal{C}([0,T], \mathcal{P}_1(\R^d))$ starting from $\mu_0$ since $b$ is bounded and Lipschitz continuous. Therefore $(\mathcal{L}(\hat{\mu}^N_.)_{n \in \mathbb{N}}$ converges to $\delta_{\mu.}$ in $\mathcal{P}( \mathcal{C}([0,T], \mathcal{P}_1(\R^d))$. The limit is deterministic and therefore $\hat{\mu}_.^N$ actually converges toward $\mu_.$ in probability. Since 
$$\sup_{N \in \mathbb{N}} \E \left[ \sup_{t \in [0,T]} \frac{1}{N} \sum_{i=1}^N | X_t^{i,N} | \right] <+\infty$$ 
we can improve the convergence from convergence in probability to convergence in $L^1$ to deduce \eqref{ToproveintheAppendixConcentration}.

\end{proof}

\paragraph{Acknowledgment} The author thanks Pierre Cardaliaguet for suggesting this problem and for fruitful discussions during the preparation of this work which is part of his PhD thesis.

\bibliographystyle{plain}
\bibliography{/Users/sam/Documents/Bibtex/DraftPaper1.bib}

\begin{thebibliography}{10}

\bibitem{Ambrosio2000}
Luigi Ambrosio, Nicola Fusco, and Diego Pallara.
\newblock {\em {Functions of Bounded Variation and Free Discontinuity
  Problems}}.
\newblock Oxford Mathematical Monographs, 2000.

\bibitem{Bonnet2019}
Beno{\^{i}}t Bonnet.
\newblock {A Pontryagin Maximum Principle in Wasserstein Spaces for Constrained
  Optimal Control Problems}.
\newblock {\em ESAIM - Control, Optimisation and Calculus of Variations},
  25:1--35, 2019.

\bibitem{Bonnet2021}
Beno{\^{i}}t Bonnet and H{\'{e}}l{\`{e}}ne Frankowska.
\newblock {Necessary Optimality Conditions for Optimal Control Problems in
  Wasserstein Spaces}.
\newblock {\em Applied Mathematics and Optimization}, pages 1--34, 2021.

\bibitem{Bouchard2009}
Bruno Bouchard, Romuald Elie, and Cyril Imbert.
\newblock {Optimal control under stochastic target constraints}.
\newblock {\em SIAM Journal on Control and Optimization}, 48(5):3501--3531,
  2009.

\bibitem{Bouchard2010}
Bruno Bouchard, Romuald Elie, and Nizar Touzi.
\newblock {Stochastic Target Problems with Controlled Loss}.
\newblock {\em SIAM Journal on Control and Optimization}, 48(5):3123--3150,
  2010.

\bibitem{Briani2018}
Ariela Briani and Pierre Cardaliaguet.
\newblock {Stable solutions in potential mean field game systems}.
\newblock {\em Nonlinear Differential Equations and Applications}, 25(1):1--26,
  2018.

\bibitem{Budhiraja2000}
Amarjit Budhiraja and Paul Dupuis.
\newblock {A variational representation for positive functionals of infinite
  dimensional Brownian motion}.
\newblock {\em Probability and Mathematical Statistics}, 20:1--25, 2000.

\bibitem{Budhiraja2012}
Amarjit Budhiraja, Paul Dupuis, and Markus Fischer.
\newblock {Large deviation properties of weakly interacting processes via weak
  convergence methods}.
\newblock {\em The Annals of Probability}, 40(1):74--102, 2012.

\bibitem{Burzoni2020}
Matteo Burzoni, Vincenzo Ignazio, Max Reppen, and H.~Mete Soner.
\newblock {Viscosity Solutions for Controlled McKean--Vlasov Jump-Diffusions}.
\newblock {\em SIAM Journal on Control and Optimization}, 58(3):1676--1699,
  2020.

\bibitem{Caines2006}
Peter~E. Caines, Minyi Huang, and Roland~P. Malham{\'{e}}.
\newblock {Large population stochastic dynamic games: closed-loop McKean-Vlasov
  systems and the Nash certainty equivalence principle}.
\newblock {\em Communications in Information and Systems}, 6(3):221--252, 2006.

\bibitem{Caines2007}
Peter~E. Caines, Minyi Huang, and Roland~P. Malham{\'{e}}.
\newblock {Large-Population Cost-Coupled LQG Problems With Nonuniform Agents :
  Individual-Mass Behavior and Decentralized epsilon-Nash Equilibria}.
\newblock {\em IEEE Trans. on Automatic Control}, 52(9):1560--1571, 2007.

\bibitem{Cardaliaguet2022}
Pierre Cardaliaguet, Samuel Daudin, Joe Jackson, and Panagiotis Souganidis.
\newblock {An algebraic convergence rate for the optimal control of
  Mckean-Vlasov dynamics}.
\newblock {\em arXiv:2203.14554}, pages 1--28, 2022.

\bibitem{Cardaliaguet2019a}
Pierre Cardaliaguet, Fran{\c{c}}ois Delarue, Jean-Michel Lasry, and
  Pierre-Louis Lions.
\newblock {The master equation and the convergence problem in mean field
  games}.
\newblock {\em Annals of Mathematics Studies}, 2019-Janua(201):1--222, 2019.

\bibitem{Cardaliaguet2016}
Pierre Cardaliaguet, Alp{\'{a}}r M{\'{e}}sz{\'{a}}ros, and Filippo
  Santambrogio.
\newblock {First order mean field games with density constraints: Pressure
  equals price}.
\newblock {\em SIAM Journal on Control and Optimization}, 54(5):2672--2709,
  2016.

\bibitem{Cardaliaguet2022a}
Pierre Cardaliaguet and Panagiotis Souganidis.
\newblock {Regularity of the value function and quantitative propagation of
  chaos for mean field control problems}.
\newblock {\em arXiv:2204.01314}, 1(1):1--26, 2022.

\bibitem{Carmona2018a}
Ren{\'{e}} Carmona and Fran{\c{c}}ois Delarue.
\newblock {\em {Probabilistic theory of mean field games with applications I.
  Mean Field FBSDEs, Control, and Games}}.
\newblock Springer Cham, 2018.

\bibitem{Carmona2018b}
Ren{\'{e}} Carmona and Fran{\c{c}}ois Delarue.
\newblock {\em {Probabilistic Theory of Mean Field Games with Applications II.
  Mean Field Games with Common Noise and Master Equations}}, volume~84.
\newblock Springer Cham, 2018.

\bibitem{Cecchin2022}
Alekos Cecchin and Fran{\c{c}}ois Delarue.
\newblock {Weak solutions to the master equation of potential mean field
  games}.
\newblock {\em arxiv 2204.04315}, pages 1--85, 2022.

\bibitem{Chow2020}
Yuk~Loong Chow, Xiang Yu, and Chao Zhou.
\newblock {On Dynamic Programming Principle for Stochastic Control Under
  Expectation Constraints}.
\newblock {\em Journal of Optimization Theory and Applications},
  185(3):803--818, 2020.

\bibitem{Conforti2021}
Giovanni Conforti, Richard Kraaij, and Daniela Tonon.
\newblock {Hamilton--Jacobi equations for controlled gradient flows: the
  comparison principle}.
\newblock {\em arxiv 2111.13258}, pages 1--43, 2021.

\bibitem{Cosso2021}
Andrea Cosso, Fausto Gozzi, Idris Kharroubi, Huy{\^{e}}n Pham, and Mauro
  Rosestolato.
\newblock {Master Bellman equation in the Wasserstein space: Uniqueness of
  viscosity solutions}.
\newblock {\em arXiv:2107.10535}, 2021.

\bibitem{Daudin2022}
Samuel Daudin.
\newblock {Optimal Control of Diffusion Processes with Terminal Constraint in
  Law}.
\newblock {\em Journal of Optimization Theory and Applications}, 195(1):1--41,
  2022.

\bibitem{Daudin2021}
Samuel Daudin.
\newblock {Optimal control of the Fokker-Planck equation under state
  constraints in the Wasserstein space}.
\newblock {\em (to appear in) Journal de Math{\'{e}}matiques Pures et
  Appliqu{\'{e}}es}, 2023.

\bibitem{Daudin2023a}
Samuel Daudin, Fran{\c{c}}ois Delarue, and Joe Jackson.
\newblock {On the Optimal Rate for the Convergence Problem in Mean Field
  Control}.
\newblock {\em arXiv:2305.08423}, 2023.

\bibitem{Dawson1987}
Donald~A. Dawson and J{\"{u}}rgen G{\"{a}}rtner.
\newblock {Large deviations from the mckean-vlasov limit for weakly interacting
  diffusions}.
\newblock {\em Stochastics}, 20(4):247--308, 1987.

\bibitem{DiMarino2016b}
Simone {Di Marino} and Alp{\'{a}}r M{\'{e}}sz{\'{a}}ros.
\newblock {Uniqueness issues for evolution equations with density constraints}.
\newblock {\em Mathematical Models and Methods in Applied Sciences},
  26(9):1761--1783, 2016.

\bibitem{Djete2020}
Mao~Fabrice Djete, Dylan Possama{\"{i}}, and Xiaolu Tan.
\newblock {McKean-Vlasov optimal control: limit theory and equivalence between
  different formulations}.
\newblock {\em Mathematics of Operations Research}, 2022.

\bibitem{Fischer2014}
Markus Fischer.
\newblock {On the form of the large deviation rate function for the empirical
  measures of weakly interacting systems}.
\newblock {\em Bernoulli}, 20(4):1765--1801, 2014.

\bibitem{Fischer2016}
Markus Fischer and Giulia Livieri.
\newblock {Continuous time mean-variance portfolio optimization through the
  mean field approach}.
\newblock {\em ESAIM - Probability and Statistics}, 20:30--44, 2016.

\bibitem{Fleming1977}
Wendell~H. Fleming.
\newblock {Exit probabilities and optimal stochastic control}.
\newblock {\em Applied Mathematics {\&} Optimization}, 4(1):329--346, 1977.

\bibitem{Fleming2006}
Wendell~H. Fleming and H.~Mete Soner.
\newblock {\em {Controlled Markov Processes and Viscosity Solutions}}.
\newblock Springer-Verlag New York, 2006.

\bibitem{Follmer1999}
Hans F{\"{o}}llmer and Peter Leukert.
\newblock {Quantile hedging}.
\newblock {\em Finance and Stochastics}, 3(3):251--273, 1999.

\bibitem{Frankowska2018}
H{\'{e}}l{\`{e}}ne Frankowska, Haisen Zhang, and Xu~Zhang.
\newblock {Stochastic optimal control problems with control and initial-final
  states constraints}.
\newblock {\em SIAM Journal on Control and Optimization}, 56(3):1823--1855,
  2018.

\bibitem{Frankowska2019}
H{\'{e}}l{\`{e}}ne Frankowska, Haisen Zhang, and Xu~Zhang.
\newblock {Necessary optimality conditions for local minimizers of stochastic
  optimal control problems with state constraints}.
\newblock {\em Transactions of the American Mathematical Society},
  372(2):1289--1331, 2019.

\bibitem{Germain2021}
Maximilien Germain, Huy{\^{e}}n Pham, and Xavier Warin.
\newblock {Rate of convergence for particle approximation of PDEs in
  Wasserstein space}.
\newblock {\em arxiv:2103.00837}, 2021.

\bibitem{Germain2022}
Maximilien Germain, Huy{\^{e}}n Pham, and Xavier Warin.
\newblock {A level-set approach to the control of state-constrained
  McKean-Vlasov equations: application to renewable energy storage and
  portfolio selection}.
\newblock {\em Numerical Algebra, Control and Optimization}, 2022.

\bibitem{Katsoulakis1994a}
Markos~A. Katsoulakis.
\newblock {Viscosity Solutions of Second Order Fully Nonlinear Elliptic
  Equations with State Constraints}.
\newblock {\em Indiana University Mathematics Journal}, 43(2):493--519, 1994.

\bibitem{Krokhmal2001}
Pavlo Krokhmal, Tanislav Uryasev, and Jonas Palmquist.
\newblock {Portfolio optimization with conditional value-at-risk objective and
  constraints}.
\newblock {\em The Journal of Risk}, 4(2):43--68, 2001.

\bibitem{Lacker2015}
Daniel Lacker.
\newblock {Mean field games via controlled martingale problems: Existence of
  Markovian equilibria}.
\newblock {\em Stochastic Processes and their Applications}, 125(7):2856--2894,
  2015.

\bibitem{Lacker2017}
Daniel Lacker.
\newblock {Limit theory for controlled mckean-vlasov dynamics}.
\newblock {\em SIAM Journal on Control and Optimization}, 55(3):1641--1672,
  2017.

\bibitem{Lasry1989}
Jean-Michel Lasry and Pierre-Louis Lions.
\newblock {Nonlinear elliptic equations with singular boundary conditions and
  stochastic control with state constraints - 1. The model problem}.
\newblock {\em Mathematische Annalen}, 283(4):583--630, 1989.

\bibitem{Lasry2006}
Jean-Michel Lasry and Pierre-Louis Lions.
\newblock {Jeux {\`{a}} champ moyen. II - Horizon fini et contr{\^{o}}le
  optimal}.
\newblock {\em Comptes Rendus Mathematique}, 343(10):679--684, 2006.

\bibitem{Lasry2006a}
Jean-Michel Lasry and Pierre-Louis Lions.
\newblock {Mean field games. I - The stationary case}.
\newblock {\em Comptes Rendus Mathematique}, 343(9):619--625, 2006.

\bibitem{Lasry2007}
Jean-Michel Lasry and Pierre-Louis Lions.
\newblock {Mean field games}.
\newblock {\em Japanese Journal of Mathematics}, 2(1):229--260, 2007.

\bibitem{Leonori2007}
Tommaso Leonori and Alessio Porretta.
\newblock {The boundary behavior of blow-up solutions related to a stochastic
  control problem with state constraint}.
\newblock {\em SIAM Journal on Mathematical Analysis}, 39(4), 2007.

\bibitem{Lions20062020}
Pierre-Louis Lions.
\newblock {Th{\'{e}}orie des jeux {\`{a}} champs moyen}.
\newblock {\em video lecture series at Coll{\`{e}}ge de France}, 2006.

\bibitem{Markowitz1952}
Harry Markowitz.
\newblock {Porfolio Selection}.
\newblock {\em The Journal of Finance}, 7(1):77--91, 1952.

\bibitem{Meszaros2015}
Alp{\'{a}}r M{\'{e}}sz{\'{a}}ros and Francisco~J. Silva.
\newblock {A variational approach to second order mean field games with density
  constraints : The stationary case}.
\newblock {\em Journal de Math{\'{e}}matiques Pures et Appliqu{\'{e}}es},
  104(6):1135--1159, 2015.

\bibitem{Meszaros2018}
Alp{\'{a}}r M{\'{e}}sz{\'{a}}ros and Francisco~J. Silva.
\newblock {On The Variational Formulation Of Some Stationary Second-Order Mean
  Field Games Systems}.
\newblock {\em SIAM Journal on Mathematical Analysis}, 50(1):1255--1277, 2018.

\bibitem{Oelschlager1984}
Karl Oelschlager.
\newblock {A Martingale Approach to the Law of Large Numbers for Weakly
  Interacting Stochastic Processes}.
\newblock {\em The Annals of Probability}, 12(2):458--479, 1984.

\bibitem{Orrieri2019}
Carlo Orrieri, Alessio Porretta, and Giuseppe Savar{\'{e}}.
\newblock {A variational approach to the mean field planning problem}.
\newblock {\em Journal of Functional Analysis}, 2019.

\bibitem{Patie2008}
P.~Patie and C.~Winter.
\newblock {First exit time probability for multidimensional diffusions: A
  PDE-based approach}.
\newblock {\em Journal of Computational and Applied Mathematics},
  222(1):42--53, 2008.

\bibitem{Pfeiffer2020}
Laurent Pfeiffer.
\newblock {Optimality conditions in variational form for non-linear constrained
  stochastic control problems}.
\newblock {\em Mathematical Control and Related Fields}, 10(3):493--526, 2020.

\bibitem{Pfeiffer2020a}
Laurent Pfeiffer, Xiaolu Tan, and Yu~Long Zhou.
\newblock {Duality and approximation of stochastic optimal control problems
  under expectation constraints}.
\newblock {\em SIAM Journal on Control and Optimization}, 59(5):3231--3260,
  2021.

\bibitem{Meszaros2016a}
Filippo Santambrogio and Alp{\'{a}}r M{\'{e}}sz{\'{a}}ros.
\newblock {Advection-Diffusion Equations With Density Constraints}.
\newblock {\em Analysis and PDE}, 9(3), 2016.

\bibitem{Seguret2021}
Adrien Seguret.
\newblock {Optimal control of a first order Fokker-Planck equation with
  reaction term and density constraints}.
\newblock {\em arXiv:2109.12836v2}.

\bibitem{Simons1998}
Stephen Simons.
\newblock {\em {Minimax and Monoticity}}.
\newblock Springer-Verlag Berlin Heidelberg, 1998.

\bibitem{Soner1986}
H.~Mete Soner.
\newblock {Optimal Control With State-Space Constraint I.}
\newblock {\em SIAM Journal on Control and Optimization}, 24(6):1110--1122,
  1986.

\end{thebibliography}

\end{document}